\documentclass[preprint,12pt]{elsarticle}

 \usepackage{epsfig}

 \usepackage{amssymb}
 \usepackage{amsthm}

 \usepackage{amsmath}

\journal{Journal of Computational Physics}

\def\ds{\displaystyle }
\def\eps{\varepsilon }

\def\u{{\bf u}}
\def\j{{\bf j}}
\def\B{{\bf B}}
\def\e{{\bf e}}
\def\p{{\bar{\psi}}}

\begin{document}

\begin{frontmatter}



\title{Reconstruction of the equilibrium of the plasma in a Tokamak 
and identification of the current density profile in real time
}


\author[nice]{J. Blum}
\author[nice]{C. Boulbe}
\author[nice]{B. Faugeras}
\address[nice]{Laboratoire J.A. Dieudonn\'e, UMR 6621, Universit\'e de Nice 
Sophia Antipolis, Parc Valrose, 06108 Nice Cedex 02, France}

\begin{abstract}
The reconstruction of the equilibrium of a plasma in a Tokamak is a free boundary
problem described by the Grad-Shafranov equation in axisymmetric
configuration. The right-hand side of this equation is a nonlinear
source, which represents the toroidal component of the plasma current
density. This paper deals with the identification of this
nonlinearity source from experimental measurements in real time. 
The proposed method is based on a fixed point algorithm, a finite
element resolution, a reduced basis method  and a least-square optimization formulation. 
This is implemented in a software called Equinox with which 
several numerical experiments are conducted to explore the identification problem. 
It is shown that the identification of the profile of the averaged current density and of the safety factor as a function 
of the poloidal flux is very robust. 

\end{abstract}


\end{frontmatter}


\section{Introduction}

In fusion experiments a magnetic field is used to confine a plasma 
in the toroidal vacuum vessel of a Tokamak \cite{Wesson:2004}. 
The magnetic field is produced by external coils 
surrounding the vacuum vessel and also by a current circulating in the plasma itself. 
The resulting magnetic field is helicoidal. 

Let us denote by ${\bf j}$ the current density in the plasma, 
by ${\bf B}$ the magnetic field and by $p$ the kinetic pressure. 
The momentum equation for the plasma is 
$$
\ds \rho \frac{d \u}{dt} +\nabla p=\j \times \B
$$
where $\u$ represents the mean velocity of particles and $\rho$ the mass density. 
At the slow resistive diffusion time scale \cite{Grad:1970} 
the term $\ds \rho \frac{d \u}{dt}$ can be neglected compared to 
$\nabla p$ and the equilibrium equation 
for the plasma simplifies to 
$$
\j \times \B = \nabla p
$$
meaning that at each instant in time the plasma is at 
equilibrium and the Lorentz force ${\bf j}\times{\bf B}$ balances 
the force $\nabla p$ due to kinetic pressure. 
Taking into account the magnetostatic Maxwell equations which are satisfied in the whole space 
(including the plasma) the equilibrium of the plasma in presence of a magnetic field is
described by
\begin{eqnarray}
\mu_0 {\bf j}&=&\nabla\times {\bf B},\label{current}\\
\nabla\cdot{\bf B}&=&0,\label{div}\\
{\bf j}\times{\bf B}&=&\nabla p,\label{equilibre}
\end{eqnarray}
where $\mu_0$ is the magnetic permeability of the vacuum. 
Ampere's theorem is expressed by Eq. (\ref{current}) 
and Eq. (\ref{div}) represents the conservation of magnetic
induction.
From the equilibrium equation (\ref{equilibre}) it is clear that 
\begin{equation*}
{\bf B}\cdot\nabla p=0 \mbox{ and }
{\bf j}\cdot\nabla p=0.
\end{equation*}
Therefore field lines and current lines lie on isobaric surfaces. 
These isosurfaces form a family of nested tori called magnetic
surfaces which enable to define the magnetic axis
and the plasma boundary. On the one hand the innermost magnetic
surface degenerates into a closed curve and is called magnetic axis 
and on the other hand the plasma boundary corresponds to the surface 
in contact with a limiter or to a magnetic separatrix (hyperbolic line with an X-point).

\begin{figure}
\begin{center}
\scalebox{0.6}{\input{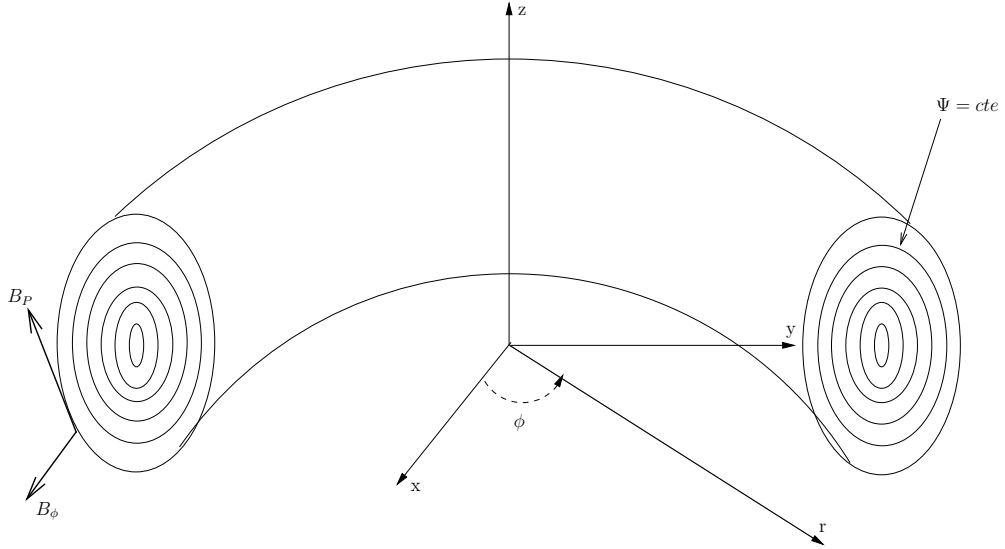}}
\caption{Toroidal geometry.}\label{fig:geomtore}
\end{center}
\end{figure}

The Grad-Shafranov equation \cite{Grad:1958,Shafranov:1958,Mercier:1974} is a rewriting of 
Eqs. (\ref{current}-\ref{equilibre}) under the axisymmetric assumption.
Consider the cylindrical coordinate system $(\e_r,\e_\phi,\e_z)$. 
The  magnetic field $\B$ is supposed to be independent of the toroidal angle $\phi$. 
Let us decompose it in a poloidal field $\B_p=B_r \e_r + B_z \e_z$ 
and a toroidal field $\B_\phi = B_\phi \e_\phi$ (see Fig. \ref{fig:geomtore}). 

Let us also introduce the poloidal flux 
$$
\psi(r,z)=\ds \frac{1}{2\pi}\int_D \B {\bf ds} = \int_0^r B_z r dr
$$
where $D$ is the disc having as circumference the circle centered 
on the $Oz$ axis and passing through a point $(r,z)$ in a poloidal section. 
From Eq. (\ref{div}) one deduces that 
$\displaystyle{{\bf B}_p}=\displaystyle{\frac{1}{r}[\nabla \psi\times {\bf e}_\phi]}$. 
Therefore $\B.\nabla \psi = 0 $ meaning that $\psi$ is a constant 
on each magnetic surface and that $p=p(\psi)$.

The same poloidal-toroidal decomposition can be applied to $\j$. From Eq. (\ref{current}) 
it is clear that 
$\nabla \cdot \j =0$. As for $\B_p$ it is shown that there exists a function $f$, called the 
diamagnetic function, such that 
$\j_p= \ds \frac{1}{r}[\nabla (\frac{f}{\mu_0}) \times \e_\phi]$. 
Since $\j.\nabla p=0$ then $\nabla f \times \nabla p =0$ and $f$ is constant on the magnetic surfaces, $f=f(\psi)$.

From Eq. (\ref{current}) one also deduces that $\B_\phi=\ds \frac{f}{r} \e_\phi$ and 
$\j_{\phi} =  (-\Delta^* \psi) \e_{\phi}$ where
\begin{equation*}
\label{eqn:delta*}
\Delta^* .=  \frac{\partial }{\partial r}(\frac{1}{\mu_0 r}  \frac{ \partial . }{\partial r}) 
+ \frac{\partial}{\partial z}(\frac{1}{\mu_0 r} \frac{ \partial .}{\partial z}).
\end{equation*}

To sum up
 
\begin{equation*}\label{decomp}
\begin{array}{lll}
\left\{
\begin{array}{lrc}
{\bf B}={\bf B}_p+{\bf B}_\phi\\
\displaystyle{{\bf B}_p}=\displaystyle{\frac{1}{r}[\nabla \psi\times {\bf e}_\phi]}\\
\displaystyle{{\bf B}_\phi}=\displaystyle{\frac{f}{r}{\bf e}_\phi}
\end{array}
\right.
&\hskip 1cm \mathrm{and} \hskip 1cm &
\left\{
\begin{array}{lrc}
{\bf j}={\bf j}_p+{\bf j}_\phi\\
\displaystyle{{\bf j}_p}=\displaystyle{\frac{1}{r}[\nabla  \frac{f}{\mu_0}  \times {\bf e}_\phi]}\\
\displaystyle{{\bf j}_\phi}=-\Delta^* \psi {\bf e}_\phi
\end{array}
\right.
\end{array}
\end{equation*}

From Eq. (\ref{equilibre}) one deduces that
$$
(\j_p+j_\phi \e_\phi) \times (\B_p + B_\phi \e_\phi) = 
\ds -\frac{1}{\mu_0 r} B_\phi \nabla f + j_\phi \ds \frac{1}{r} \nabla \psi = \nabla p
$$
and since 
$$
\nabla p= p'(\psi) \nabla \psi \ \mathrm{and}\  \nabla f= f'(\psi) \nabla \psi
$$
the Grad-Shafranov equation valid in the plasma reads
\begin{equation}
\label{eqn:gradshaf}
-\Delta^* \psi = r p'(\psi) + \frac{1}{\mu_0 r}(ff')(\psi)
\end{equation}

Thus under the axisymmetric assumption, the three dimensional equilibrium
Eqs. (\ref{current} - \ref{equilibre}) reduce to a two dimensional 
non linear problem. Note that the right-hand side of
Eq. (\ref{eqn:gradshaf}) represents the toroidal component $j_{\phi}$ of the 
current density in the plasma which is determined by the unknown functions $p'$ and $ff'$. 
In the vacuum there is no current and the poloidal flux satisfies 
$$
-\Delta^* \psi = 0
$$

In this paper, we are interested in the numerical reconstruction of the
equilibrium i.e of the poloidal flux $\psi$ and in the identification of the unknown 
plasma current density \cite{Lao:1990,Blum:1990,Blum:1997}. 
In a control perspective this reconstruction has to be achieved in real time 
from experimental measurements. 
The main difficulty consists in identifying the functions $p'$ and $ff'$ 
in the non linear right-hand side source term in Eq. (\ref{eqn:gradshaf}). 
An iterative strategy involving a finite element method for the resolution of 
the direct problem and a least square optimisation procedure for the 
identification of the non linearity using a decomposition basis is proposed.

Let us give a brief historical background of this problem of the reconstruction 
of the plasma current density from experimental measurements. 
In large aspect ratio Tokamaks with circular cross-sections, it was established in 
\cite{Shafranov:1971, Zakharov:1973} 
that the quantities that can be identified from magnetic measurements 
are the total plasma current $I_p$ 
and a sum involving the poloidal beta and the internal inductance: $\beta_p+l_i/2$ (see Appendix \ref{appendixC}). 
A large number of papers proved the possibility of separating $\beta_p$ from $l_i$ 
as soon as the plasma is no longer circular with high-aspect ratio 
\cite{Luxon:1982, Swain:1982, Lao:1985, Blum:1989}. 
The fact of adding supplementary experimental diagnostics, 
such as line integrated electronic density and Faraday rotation measurements, 
has considerably improved the identification of the current density profile 
\cite{Hofmann:1988, Lao:1990, Blum:1990}. 
The knowledge of the flux lines (from density or temperature measurements) enables 
in principle \cite{Christiansen:1982} to determine fully the two functions $p'$ and $ff'$ in the toroidal 
plasma current density, except in a particular case pointed out by \cite{Braams:1991} 
and studied by \cite{Bishop:1986} 
and referred to as minimum-B equilibria. 
The difficulty in the reconstruction of the current profile, 
especially when only magnetic measurements are used, 
has been pointed out in 
\cite{Pustovitov:2001} 
and is inherent to the ill-posedness of this inverse problem. 
The theory of variances in equilibrium reconstruction 
\cite{Zakharov:2008} 
enables to determine by statistical methods what kind of plasma functions 
can be reconstructed in a robust way. 
The equilibrium reconstruction problem in the case of anisotropic 
pressure is treated in \cite{Zwingmann:2001}. 

A certain number of mathematical results on the identifiability of the right-hand-side of the Grad-Shafranov equation 
from Cauchy boundary conditions on the plasma frontier exist and seem unknown from the physical community. 
They are first dealing with the cylindrical case where the equilibrium equation becomes
$-\Delta \psi=p'(\psi)$
and where only one non-linearity has to be identified. 
It is clear that, if the plasma boundary is circular, then the magnetic field is constant 
on the plasma boundary and there is an infinity of non-linearities giving this value and 
the only information coming from the poloidal field on the plasma boundary is the total plasma current. 
In \cite{Beretta:1991} it was proved that if $p'$ is a real-analytic function, then in a domain with a corner there is only 
one non-linearity $p'$ corresponding to a given poloidal field on the plasma boundary.  
Some angles in the proof were excluded but in \cite{Beretta:1995} the proof was given for corners with arbitrary angles 
(including the  $90$ degrees X-point case). 
Curiously the case where the plasma boundary is smooth is mathematically more difficult and it has been proved 
in \cite{Vogelius:1994} that, if the plasma is non-circular and if $p'$ is affine in terms of $\psi$ 
then there exists at most a finite number of affine functions corresponding to the Cauchy boundary conditions. 
The link with the Schiffer and Pompeiu conjectures which is clearly pointed out in this paper is particularly interesting. 
In \cite{Demidov:2009} results of unicity for a class of affine functions or for exponential functions are given for special 
smooth boundaries and results of non-unicity for doublet-type configurations. 
Finally in \cite{Beretta:1992} identifiability results are given for the full Grad-Shafranov equation in a domain with a corner, 
with some exceptions for the angle. Of course, in spite of all these identifiability results, 
the ill-posedness of the reconstruction of the non-linearities from the Cauchy boundary measurements 
remains and has to be tackled very cautiously.

Section 2 is devoted to the statement of the mathematical problem and to the description 
of the experimental measurements avalaible. The proposed algorithm is described in Section 3. 
This methodology has been implemented in a software called Equinox 
and numerical results using synthetic and real measurements are presented in Section 4.

\section{Setting of the direct and inverse problems}

\subsection{Experimental measurements}

Although the unknown functions $p'(\psi)$ and $(ff')(\psi)$ cannot be directly measured in a Tokamak several
measurements are available:

\begin{itemize}
\item Magnetic measurements: they represent the basic information on which any equilibrium reconstruction relies. 
Flux loops provide measurements of $\psi$ and 
magnetic probes provide measurements of the poloidal field ${\bf B}_p$ 
at several points around the vacuum vessel.
Let $\Omega$ be the domain representing the vacuum vessel and $\partial \Omega$ its boundary. 
In what follows we assume that we are able to obtain the Dirichlet boundary conditions $\psi=g_D$ 
and the Neumann boundary conditions $\displaystyle{\frac{1}{r}\frac{\partial\psi}{\partial n}}=g_N$ at any points of the contour $\partial \Omega$ thanks to a preprocessing 
of the magnetic measurements. 
This preprocessing can either be a simple interpolation between real measurements 
or be the result of some boundary reconstruction algorithm which computes $\psi$ outside the plasma 
satisfying $\Delta^* \psi =0$ under the constraint of the measurements \cite{OBrien:1993, Saint-Laurent:2001, Sartori:2003}. 
\end{itemize}

A second set of measurements which can be used 
as a complement to magnetic measurements are internal measurements:
\begin{itemize}
\item Interferometric measurements: they give the values of the integrals along a family of chords $C_i$ of the 
electronic density $n_e(\psi)$ which is approximately constant on each flux line
$
\displaystyle{\int_{C_i} n_e(\psi) \,dl=\gamma_i}.
$

\item Polarimetric measurements: they give the value of the integrals 
$$
\int_{C_i} \frac{n_e(\psi)}{r}\frac{\partial\psi}{\partial n}dl=\alpha_i.
$$
$\displaystyle{\frac{\partial\psi}{\partial n}}$ is the
normal derivative of $\psi$ along the chord $C_i$.
\end{itemize}

Even when using magnetic measurements only for the equilibrium reconstruction 
the numerical algorithm presented in this paper also uses:  
\begin{itemize}
\item Current measurement: it gives the value of the total plasma current $I_p$ defined by 
$$
I_p=\int_{\Omega_p} j_\phi dx.
$$
Ampere's theorem shows that this quantity can be deduced from magnetic measurements.
 
\item Toroidal field measurement: it gives the value $B_0$ of the toroidal component of the field in the vacuum at
the point $(R_0,0)$ where $R_0$ is the major radius of the Tokamak. 
This is used for the integration of $ff'$ into $f$ and for the computation of 
the safety factor $q$ (see Appendix \ref{appendixA}). 
\end{itemize}


\subsection{Direct problem}
 
The equilibrium of a plasma in a Tokamak is a free boundary problem. 
The plasma boundary is determined either as being the last flux line in a
limiter $L$ or as being a magnetic separatrix with an X-point (hyperbolic point). 
The region $\Omega_p\subset \Omega$ containing the plasma
is defined by
$$
\Omega_p=\{{\bf x}\in\Omega,\;\psi({\bf x})\geq \psi_b\}
$$
where $\psi_b=\max_L \psi$ in the limiter configuration or
$\psi_b=\psi(X)$ when an X-point exists.

In the vacuum region, the right-hand side of Eq. (\ref{eqn:gradshaf}) vanishes
and the equilibrium equation reads
$$
\Delta^*\psi=0 \mbox{ in }\Omega\setminus\Omega_p
$$
Let us introduce the normalized flux 
$\bar{\psi}= \displaystyle \frac{\psi -\psi_a}{\psi_b - \psi_a} \in [0,1]$ in $\Omega_p$ with\\ 
$\psi_a=\max_{\Omega_p} \psi$, 
$A(\bar{\psi})= \displaystyle \frac{R_0}{\lambda}p'({\psi})$ and 
$B(\bar{\psi})=\displaystyle \frac{1}{\lambda\mu_0 R_0}(ff')({\psi})$. 
This is introduced so that the non dimensional and unknown 
functions $A$ and $B$ are defined and identified on the
fixed interval $[0,1]$. 
Imposing Dirichlet boundary conditions the final equilibrium equation is expressed as the 
boundary value problem:
\begin{equation}
\label{eqn:final}
\left \lbrace
\begin{array}{rcl}
-\Delta^* \psi &= &\lambda[\displaystyle\frac{r}{R_0} A(\bar{\psi}) +  \displaystyle 
\frac{R_0}{r} B(\bar{\psi})] \chi_{\Omega_p}\quad \mathrm{in}\ \Omega \\[10pt]
\psi&=& g_D\quad  \mathrm{on}\  \partial \Omega 
\end{array} 
\right.
\end{equation}
The free boundary aspect of the problem reduces to the particular non linearity appearing through 
$\chi_{\Omega_p}$ the characteristic function of $\Omega_p$. 
The parameter $\lambda$ is a scaling factor used to ensure that the given 
total current value $I_p$ is satisfied
\begin{equation}\label{lambda}
I_p=\lambda \int_{\Omega_p}[\displaystyle\frac{r}{R_0} A(\bar{\psi}) +  \displaystyle 
\frac{R_0}{r} B(\bar{\psi})]  dx.
\end{equation}

\subsection{Inverse problem}

The inverse problem consists in the identification of functions $A$ and $B$ from the 
measurements available. It is formulated as a least-square minimization problem 

\begin{equation}
\left \lbrace
\begin{array}{l}
\mathrm{Find}\ A^*,\ B^*,\ n_e^*\ \mathrm{such}\ \mathrm{that}: \\[10pt]
J(A^*,B^*,n_{e}^*)=\inf J(A,B,n_{e}). 
\end{array}
\right.
\end{equation}

If magnetic measurements only are used the formulation only needs the $A$ and $B$ 
variables and the $J_1$ and $J_2$ terms 
in Eq. (\ref{eqn:costfunction}) below are not needed. 
When polarimetric and interferometric measurements are used, the electronic density $n_e(\bar\psi)$ 
also has to be identified even if it does not appear in Eq. (\ref{eqn:final}).
The cost function $J$ is defined by
\begin{equation}
\label{eqn:costfunction}
J(A,B,n_e)=J_0 + J_1 + J_2 + J_\eps
\end{equation}
$J_0$ describes the misfit between computed and measured tangential component of $\bf{B}_p$
$$
J_0=  \displaystyle \ds \frac{1}{2} \sum_{k=1}^N \displaystyle (w_k)^2( \displaystyle \frac{1}{r} \frac{\partial \psi}{\partial n}(M_k) - g_N(M_k))^2 
$$
where $N$ is the number of points $M_k$ of the boundary $\partial \Omega$ where the magnetic measurements are given. 

$$ 
J_1=  \displaystyle \ds \frac{1}{2} \sum_{k=1}^{N_c} \displaystyle (w^{polar}_k)^2(\displaystyle \int_{C_k} 
\frac{n_e(\p)}{r} \frac{\partial \psi}{\partial n}dl - \alpha_k)^2
$$
and 
$$
J_2= \displaystyle \ds \frac{1}{2} \sum_{k=1}^{N_c} \displaystyle (w^{inter}_k)^2 \displaystyle (\displaystyle \int_{C_k} n_e(\p)dl - \gamma_k)^2
$$

$N_c$ is the number of chords over which interferometry and polarimetry measurements are given. 
The weights $w$ give the relative importance of the different measurements used. 
The influence of the choice of the weights on the results of the identification was extensively studied in 
\cite{Blum:1990}.
As a consequence of the ill-posedness of the identification of $A$, $B$ and $n_e$, 
a Tikhonov regularization term $J_\eps$ is introduced \cite{Tikhonov:1977} where
$$
J_\eps = 
\ds \frac{\eps_A}{2} \displaystyle \int_0^1 [A''(x)]^2 dx 
+
\ds \frac{\eps_B}{2} \displaystyle \int_0^1 [B''(x)]^2 dx 
+
\ds \frac{\eps_{n_e}}{2} \displaystyle \int_0^1 [n_e''(x)]^2 dx 
$$
and $\eps_A$, $\eps_B$ and $\eps_{n_e}$ are the regularization parameters.

The values of the different weights and parameters introduced in the cost 
function are discussed in Section \ref{sec:results}. 

It should be noticed here that magnetic measurements provide 
Dirichlet and Neumann boundary conditions. The choice was made to use 
the Dirichlet boundary conditions in the resolution of direct problem 
and to include the Neumann boundary conditions in the cost function formulated 
to solve the inverse problem. This is arbitrary and another solution could have been chosen.

\section{Algorithm and numerical resolution}
 
\subsection{Overview of the algorithm}
The aim of the method is to reconstruct the equilibrium and the toroidal current
density in real time. 
At each time step determined by the availability of new measurements during a discharge, the algorithm 
consists in constructing a sequence $(\psi^n, \Omega_p^n, A^n, B^n,\lambda^n)$ 
converging to the solution vector $(\psi, \Omega_p, A, B,\lambda)$. 
The unknown function $n_e$ may be added too if interferometry and polarimetry measurements are used. 
The sequence is obtained through the following iterative loop:

\begin{itemize}
\item Starting guess: $\psi^0$, $\Omega_p^0$, $A^0$, $B^0$ and $\lambda^0$ known from the previous time step solution. 

\item Step 1 - Optimisation step: compute $\lambda^{n+1}$ satisfying (\ref{lambda})
$$
\lambda^{n+1}=I_p /  \int_{\Omega_p^n}[\displaystyle\frac{r}{R_0} A^n(\bar{\psi}^n) 
+  \displaystyle \frac{R_0}{r} B^n(\bar{\psi}^n)]  dx
$$ 
then compute $A^{n+1}(\bar{\psi}^n)$ and $B^{n+1}(\bar{\psi}^n)$ using 
the least square procedure detailed in Section \ref{sec:detailedalgo}.

\item Step 2 - Direct problem step: compute $\psi^{n+1}$ solution to
\begin{equation}\label{GSA}
\left \lbrace
\begin{array}{rcl}
-\Delta^* \psi^{n+1} &= &\displaystyle{\lambda^{n+1} [\frac{r}{R_0} A^{n+1}(\bar{\psi}^{n})} +  \displaystyle 
\frac{R_0}{r} B^{n+1}(\bar{\psi}^{n})] \chi_{\Omega_p^n}\quad \mathrm{in}\ \Omega \\[10pt]
\psi^{n+1}&=& g_D\quad  \mathrm{on}\  \partial \Omega. 
\end{array} 
\right.
\end{equation}
and the new plasma domain $\Omega_p^{n+1}$. 

\item $n:=n+1$. If the process has not converged return to Step 1 
else $(\psi, \Omega_p, A, B, \lambda)=(\psi^n, \Omega_p^{n}, A^n, B^n, \lambda^{n})$. 
The process is supposed to have converged when the relative 
residu $\displaystyle{\frac{||\psi^{n+1}-\psi^{n}||}{||\psi^{n}||}}$ is small enough.
\end{itemize}

At each iteration of the algorithm, an inverse problem corresponding to the optimization step and an 
approximated direct Grad-Shafranov problem have to be solved successively. 
In Eq. (\ref{GSA}), $\bar{\psi}^n$ is known and since the right-hand side does not depend on
 $\psi^{n+1}$ the boundary value problem (\ref{GSA}) is linear.

In the next section the numerical methods used to solve the two problems corresponding to step 1 and step 2 are detailed.

\subsection{Numerical resolution}
\label{sec:numerical-resolution}

\subsubsection{The finite element method for the direct problem}
The resolution of the direct problem is based on a classical $P^1$ finite element
method \cite{Ciarlet:1980}. 
Let us consider the family of triangulation ${\tau}_h$ of
$\Omega$, and $V_h$ the finite dimensional subspace of $H^1(\Omega)$
defined by
$$
V_h=\{v_h\in H^1(\Omega), v_{h|T}\in P^1(T),\,\forall T\in {\tau}_h\}.
$$
and introduce $V_h^0=V_h\cap H^1_0(\Omega)$. The discrete variational
formulation of the boundary value problem (\ref{GSA}) reads
\begin{equation}
\label{eqn:var}
\left \lbrace
\begin{array}{l}
\mathrm{Find}\ \psi_h\in V_h \mbox{ with }
\psi_h=g_D\mbox{ on }\partial\Omega \mbox{ such that }\\[10pt]
\displaystyle \forall v_h \in V_h^0, \int_\Omega \displaystyle \frac{1}{\mu_0 r}\nabla \psi_{h} \cdot \nabla v_h dx = 
\int_{\Omega_p}\lambda [\displaystyle \frac{r}{R_0} A(\bar{\psi^*})+
  \displaystyle \frac{R_0}{r}B(\bar{\psi^*})] v_h
dx \\[10pt]
\end{array}
\right .
\end{equation}
where $\psi^*$ represents the known value of $\psi$ at the previous iteration. 
Numerically the Dirichlet boundary conditions are imposed using the method
consisting in computing the stiffness matrix $\hat{K}$ of the Neumann problem and
modifying it. 
Consider $(v_i)$ a basis of $V_h$ then 
$\hat{K}_{ij}=\ds \int_\Omega \ds \frac{1}{\mu_0 r} \nabla v_i \nabla v_j dx$
The modifications consist in replacing the rows corresponding
to each boundary node setting $1$ on
the diagonal terms and $0$ elsewhere. 
At each iteration only the right-hand side of the linear system in which the Dirichlet boundary conditions appear has to be modified. 
The linear system corresponding to Eq. (\ref{eqn:var}) can be written in the form

\begin{equation}\label{linear}
K.\Psi=y+g
\end{equation}

where $K$ is the $n\times n$ modified stiffness matrix, $\Psi$ is the unknown 
vector of size $n$ (the number of nodes 
of the finite elements mesh), $y$ is the vector associated 
with the modified right-hand side of 
Eq. (\ref{eqn:var}) 
and $g$ is the vector corresponding to the Dirichlet boundary conditions.

The matrix $K$ is sparse and let $LU$ be its decomposition. 
The inverse matrix $K^{-1}$ is not sparse. 
The linear system (\ref{linear}) is inverted using the $LU$ decomposition 
since it is computationally cheaper than using 
the full inverse matrix $K^{-1}$ which is nevertheless needed for the optimization 
step of the algorithm in Eq. (\ref{eqn:normal}) below.

The vector $y$ depends on functions $A$ and $B$ which are determined in
the optimization step. Functions $A$, $B$ and $n_e$ are decomposed on a finite dimensional basis $(\Phi_i)_{i=1,...,m}$ of 
functions defined on $[0,1]$ 
$$
A(x)=\sum_i^m a_i\Phi_i(x),\ B(x)=\sum_i^m b_i\Phi_i(x)\ \mathrm{and}\  
n_e(x)=\sum_i^m c_i\Phi_i(x).
$$
The vector $y$ reads
\begin{equation}\label{linearRhs}
y=Y(\bar{\psi^*})u
\end{equation}
where $u=(a_1,...,a_m,b_1,...,b_m)\in\mathbb{R}^{2m}$ is the vector of the components of functions $A$ and $B$ in the basis $(\Phi_i)$. 
The matrix $Y$ of size $n\times 2m$ is defined as follows. 
Each row $i$ of $Y$ is decomposed as
$$
Y_{ij}(\bar\Psi^*)=\left\{
\begin{array}{llc}
\displaystyle{\int_{\Omega_p} \lambda \displaystyle\frac{r}{R_0}\Phi_j(\bar{\psi^*}) v_i dx} &\mbox{ if }1\leq j\leq m\\
\displaystyle{\int_{\Omega_p} \lambda \displaystyle\frac{R_0}{r}\Phi_{j-m}(\bar{\psi^*}) v_i dx} &\mbox{ if }m+1\leq j\leq 2m.
\end{array}
\right.
$$

\subsubsection{Detailed numerical algorithm}
\label{sec:detailedalgo}
One equilibrium computation corresponds to one instant in time during a pulse. 
The quasi-static approximation consists in considering that at each instant the Grad-Shafranov equation is satisfied. 
During a pulse successive equilibrium configurations are computed with a time resolution 
$\Delta t$ corresponding to the acquisition time of measurements:

\begin{itemize}
\item Initialization before the discharge: the modified stiffness matrix $K$, its $LU$ decomposition as well its inverse $K^{-1}$ are computed 
once for all and stored.

\item Consider that the equilibrium at time $t-\Delta t$ is known and that a new set of measurements is acquired at time $t$.

\item Computation of the new equilibrium at time $t$ through the iterative loop briefly described in the previous Section and detailed below:
\end{itemize}

The equilibrium from the previous time step is used as a first guess in the iterative loop.

\paragraph{Step 1 - Optimization step}
During the optimisation step, $n_e$ is first estimated from interferometric measurements and $A$ and $B$ are computed in a second time.
\begin{itemize}
\item Compute the electronic density $n_e$ based on the equilibrium of the previous iteration $\p^*$ using a least square
formulation for the minimun of $J_2$ with Tikhonov regularization and solving the associated normal equation:
The flux $\p^*$ is given.
$$
n_e(x)=\ds \sum_{j=1}^{m} v_j \phi_j(x)
$$
The interferometric measurements for $i=1\ ...\ n_c$ are 
$$
\gamma_i \approx \ds \int_{C_i} n_e(\p^*) dl = \sum_j v_j \int_{C_i} \phi_j(\p^*) dl =\sum_j v_j B_{ij}
$$
The cost functional reads
$$
\begin{array}{lll}
J(v)&=&\ds \frac{1}{2}\sum_i (w_i^{inter})^2 (\sum_j B_{ij} v_j - \gamma_i)^2 + \frac{\eps}{2} v^T \Lambda v \\
    &=&\ds \frac{1}{2}||D^{1/2}(Bv - \gamma)||^2 + \frac{\eps}{2} v^T \Lambda v
\end{array}
$$
where $D^{1/2}=diag(w_i^{inter})$ and the regularization matrix $\Lambda$ 
is defined by
$$
\Lambda_{ij}=\int_0^1 \Phi"_i(x)\Phi"_j(x)dx
$$
and $\Phi^"_i$ is the second order derivative of the basis function $\Phi_i$.

It is minimized solving the associated normal equation
\begin{equation}
\label{eqn:nenormal}
(\alpha^2(D^{1/2}B)^T(D^{1/2}B) + \hat{\eps} \Lambda )\hat{v} = \alpha(D^{1/2}B)^T D^{1/2} \gamma
\end{equation}
Since $n_e \approx 10^{19} m^{-3}$ an adimensionalizing parameter $\alpha = 10^{19} m^{-3}$, 
such that $v=\alpha \hat{v}$,
is introduced in order to precondition the linear system which is inverted using LU decomposition, 
as well as a reasonable prescribed value for the non dimensional regularization parameter 
$\hat{\eps}=\alpha^2 \eps$. 


\item Compute $\lambda^{n+1}$ satisfying Eq. (\ref{lambda}). In the right-hand side $y$, $\lambda$ appears in the product $\lambda u$. 
In order to avoid any divergence issue due to the non uniqueness of $\lambda$ (for all $\alpha$, $\lambda u = (\lambda \alpha) \ds (\frac{u}{\alpha})$) 
the degrees of freedom (dofs) $u$ are scaled by $m=\max(|a_i|)$, $u$ is replaced by $\ds \frac{1}{m}u$ and $\lambda$ by $m\lambda$.

\item Compute $A$ and $B$. In order to approximate $A$ and $B$, suppose $n_e$ is known and consider the discrete approximated inverse problem
\begin{equation}
\label{eqn:optim}
\left \lbrace
\begin{array}{l}
\mathrm{Find}\ u\ \mathrm{minimizing}:\\
J(u)=\displaystyle\frac{1}{2}
\|  C({\psi^*})\Psi - d \|^2_D 
+\displaystyle \frac{\eps}{2}u^{T} \Lambda u
\end{array}
\right.
\end{equation}
where $C(\psi^*)$ is the observation operator and $d$ the vector of experimental measurements. 
The first term in $J$ is the discrete version of $J_0+J_1$. 
The second one corresponds to the first two terms of the
Tikhonov regularization $J_\eps$ with $\eps_A = \eps_B=\eps$ which will always be assumed 
in order for functions $A$ and $B$ to play a symmetric role.

Let us denote by $l$ the number of measurements available ($l=N+N_c$ if magnetic and polarimetric measurements are used) 
and by $D$ the diagonal matrix made of the weights $w_k$ and $w_k^{polar}$, the norm $\|.\|_D$ is defined by 
$\forall {\bf x}\in\mathbb{R}^l\;\; \|{\bf x}\|_D^2=(D {\bf x},{\bf x})=(D^{1/2}{\bf x},D^{1/2}{\bf x})$

$C(\psi^*)$ is a sparse matrix of size $l\times n$ and can be viewed as a vector
composed of two blocks $C_0$ of size $N\times n$ and independent of $\psi^*$ 
and $C_1(\psi^*)$ of size $N_c\times n$ corresponding respectively to $J_0$ and $J_1$. 
That is to say that multiplication of the kth row of $C_0$ by $\psi$ gives the kth Neumann boundary condition approximation 
$$(C_0)_k \Psi \approx  \ds (\frac{1}{r}\frac{\partial \psi}{\partial n})(M_k).$$  
The block $C_1(\psi^*)$ depends on $\psi^*$ 
through the $n_e(\psi^*)$ function. The multiplication of the kth row of $C_1(\psi^*)$ by $\Psi$ gives the kth polarimetric 
measurements approximation 
$$ (C_1(\psi^*))_k \Psi \approx \ds \int_{c_k} \frac{n_e(\psi^*)}{r}\frac{\partial\psi}{\partial n}dl .$$

The matrix $\Lambda$ is of size $2m\times 2m$ and is block diagonal composed of two blocks $\Lambda_1$ and $\Lambda_2$ of size $m\times m$, with
$$
(\Lambda_1)_{ij}=(\Lambda_2)_{ij}=\int_0^1 \Phi"_i(x)\Phi"_j(x)dx
$$
Using Eqs .(\ref{linear} - \ref{linearRhs}) problem (\ref{eqn:optim}) becomes
$$
\begin{array}{lll}
J(u)&=&\displaystyle\frac{1}{2}\|  C({\psi^*})\Psi - d \|_D^2 + \displaystyle\frac{\eps}{2}u^{T} \Lambda u \\[10pt]
&=&\displaystyle\frac{1}{2}\|  C({\psi^*}) {K}^{-1} Y(\bar{\psi^*}) u +( C(\psi^*){K}^{-1} g - d ) \|^2_D 
+\displaystyle\frac{\eps}{2} u^{T} \Lambda u \\[10pt]
&=&\displaystyle\frac{1}{2}\|  E u- f \|^2_D + \displaystyle\frac{\eps}{2}u^{T} \Lambda u 
\end{array}
$$
where $E= C({\psi^*}) {K}^{-1} Y(\bar{\psi^*})$ and $f= -C({\psi^*}){K}^{-1} h + d$.
Setting $\tilde{E}=D^{1/2}E$, problem (\ref{eqn:optim}) reduces to solve the normal equation
\begin{equation}
\label{eqn:normal}
(\tilde{E}^T \tilde{E}+\eps\Lambda)u=\tilde{E}^T f
\end{equation}
whose solution is denoted by $u^*$.

\end{itemize}

\paragraph{Step 2 - Direct problem step.}
Update the dofs $u$ and update the flux $\psi$ by solving the linear system
\begin{equation}
\label{eqn:dir}
K\psi= Y(\p^*)u^*+g
\end{equation}
using the $LU$ decomposition of matrix $K$. Update $\Omega_p$ possibly computing 
the position of the X-point if the plasma is not in a limiter configuration.

~\\
Finally it should be noticed that this algorithm is particularly well adapted to real-time applications. Indeed 
during the computations the expensive operations are the updates of matrices $C$ and $Y$ 
as well as the computation of products $CK^{-1}$ and $CK^{-1}Y$ which appear in Eq. (\ref{eqn:normal}). In order to reduce 
computation time the $K^{-1}$ matrix is precomputed and only the $\psi$-dependent part of $C$ is 
dealt with.
The resolution of the direct problem, Eq. (\ref{eqn:dir}), is cheap since the $LU$ decomposition 
of the $K$ matrix is also precomputed.

\section{Numerical results}
\label{sec:results}

\subsection{Twin experiment with noise free magnetic measurements}

\begin{figure}[!h]
\begin{center}
  \includegraphics[width=8cm]{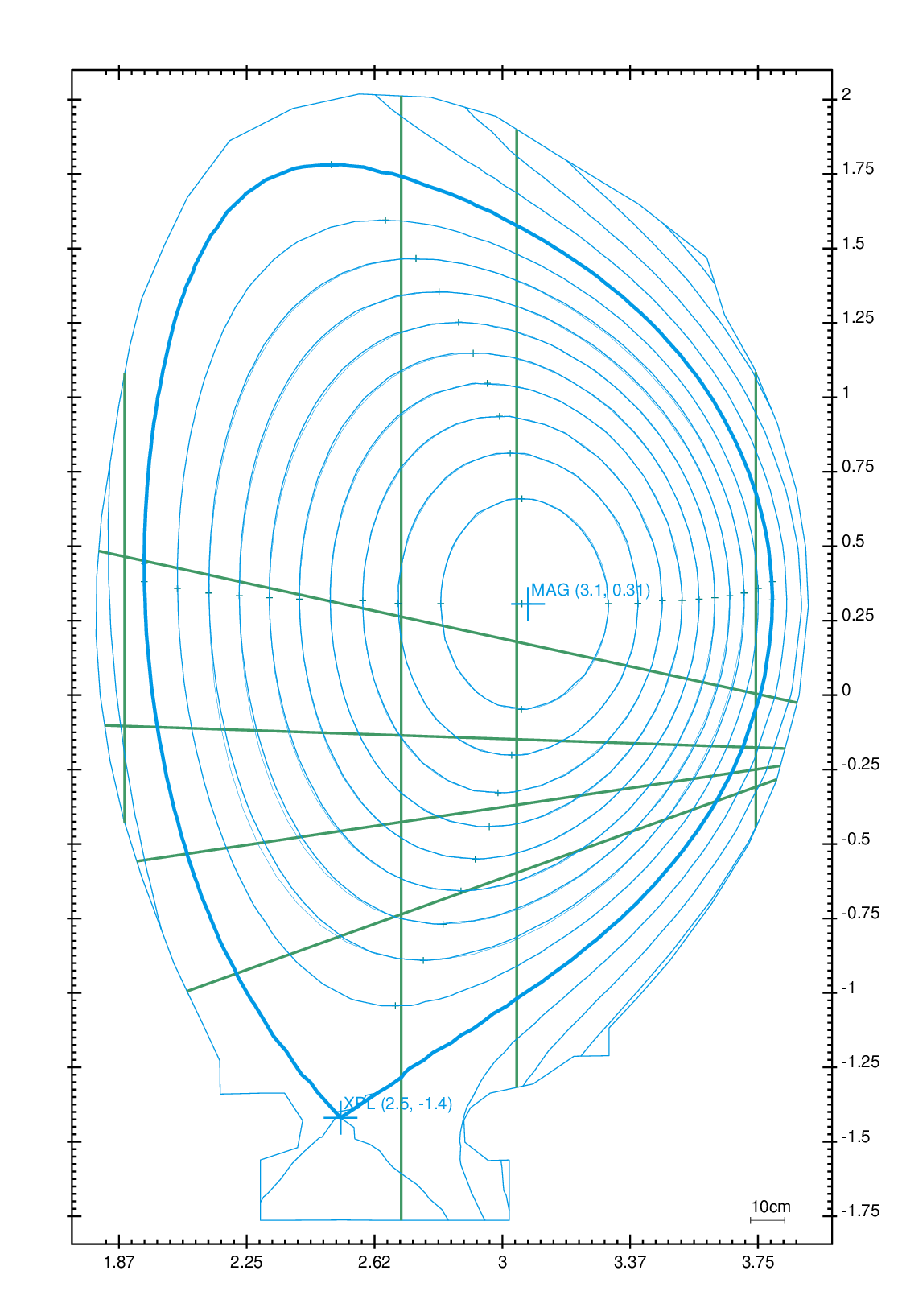} 
 \end{center}
 \caption{An equilibrium configuration for the tokamak JET 
from which twin experiments are performed. The domain $\Omega$ and its boundary $\Gamma$ (blue line) are shown. 
Isoflux are plotted from $\p=0$ (magnetic axis) to $\p=1$ 
(plasma boundary defined by the existence of an X-point at point $r=2.5$ and $z=-1.4$ m) 
by step of $\Delta \p= 0.1$ 
Interferometry and polarimetry chords appear in green.   
\label{fig:fluxmap-epsilon}}
\end{figure}

In this section we assume that the poloidal flux 
corresponding to an equilibrium configuration $\psi$ is given 
on the boundary $\Gamma$ . 
These Dirichlet boundary conditions can either 
be real measurements or can be the output from some 
equilibrium simulation code. In a first step we also assume to know 
functions $p'$ and $ff'$ (or $A$ and $B$). In what follows these reference 
functions are given point by point. 
It is then possible to run a direct simulation to compute $\psi$ on $\Omega$ 
(see Fig. \ref{fig:fluxmap-epsilon}) 
and thus 
$\ds \frac{1}{r}\frac{\partial \psi}{\partial n}$ on $\Gamma$ 
which can then be used as measurements in an inverse problem resolution. 

In this first experiment the magnetic measurements are free of noise. 
The identification algorithm is initialized using the first guess functions are 
$A(x)=B(x)=1-x$ and $\lambda=1$. The poloidal flux $\psi$ is initially a constant on $\Omega$. 
The weights in the misfit part of the cost function $J_0$ 
related to magnetic measurements are defined by 
$w_k=\ds \frac{1}{\sqrt{N} \sigma}$. 
Since the error on magnetic measurements are of about one percent we define 
$\sigma=0.01 B_m$ where $B_m$ is a mean magnetic field value which thanks 
to Ampere's theorem can be defined as $B_m = \ds \frac{\mu_0 I_p}{|\Gamma|}$.

The functions $A$ and $B$ are decomposed in a function 
basis defined on the interval $[0,1]$. Several basis have been tested 
(piecewise affine functions, polynomials, B-splines and wavelets) 
in order to verify that the result of the identification does not depend on the 
decomposition basis.  This is the case as long as the dimension of the basis 
is large enough. In the remaining part of this paper each function is decomposed 
in the same basis of 8 B-splines \cite{DeBoor:1978}. 
The boundary condition $A(1)=B(1)=0$ is imposed.

The computations are carried out for several values 
of the regularization parameters $\eps$ ranging from $10^{-10}$ to $1$. 
We are interested in the ability of the method to recover 
functions $A$ and $B$ and thus the current density profile averaged 
over the magnetic surfaces (see Appendix \ref{appendixA}):  
$$
 R_0<\ds \frac{j(r,\p)}{r}>= \lambda A(\p) + \lambda R_0^2 <\frac{1}{r^2}>B(\p)
$$
and the safety factor $q$ (see Appendix \ref{appendixB}).

As can be seen from Fig. \ref{fig:twinAB} 
the optimal choice for $\eps$ is of about $10^{-5}$ 
for which functions $A$ and $B$ are well recovered. 
For smaller values some oscillations appear because the regularization 
is not strong enough and on the contrary greater values lead 
to less precision in the recovery of the unknown functions since 
regularization is too strong. 
In the second column the relative errors on the identified functions are plotted.

\begin{figure}[!h]
\begin{center}
\begin{tabular}{ll}
  \includegraphics[width=6cm]{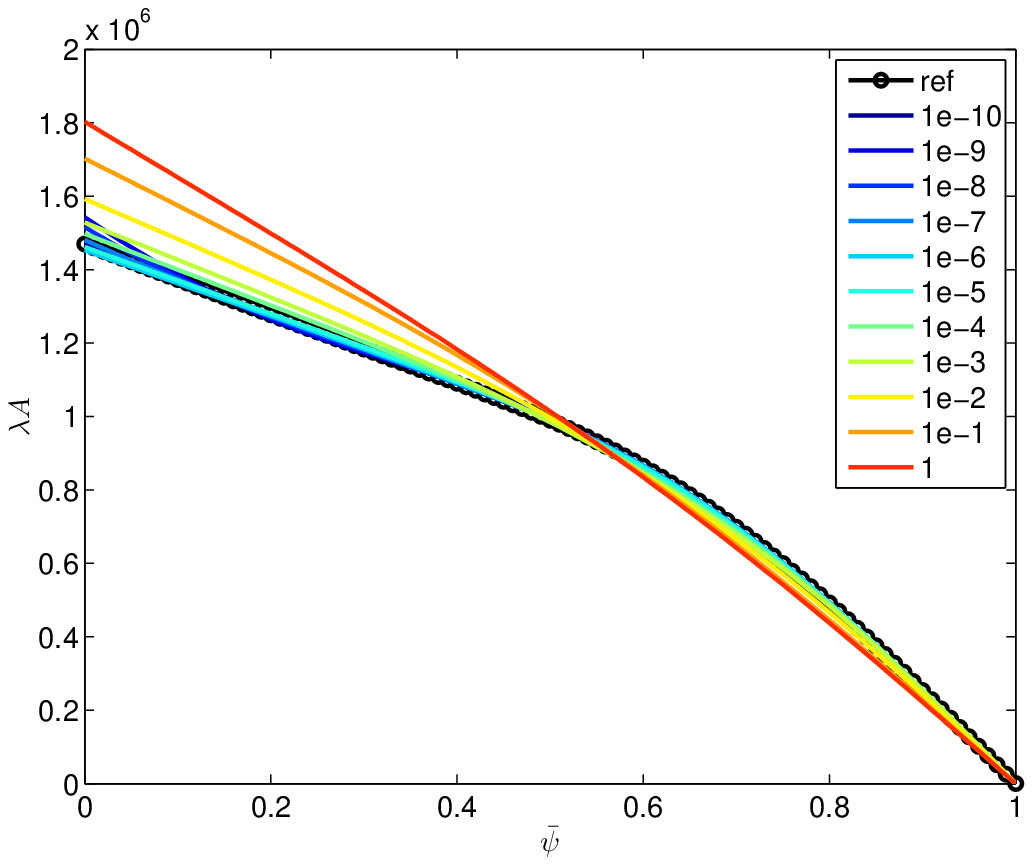} &  
  \includegraphics[width=6cm]{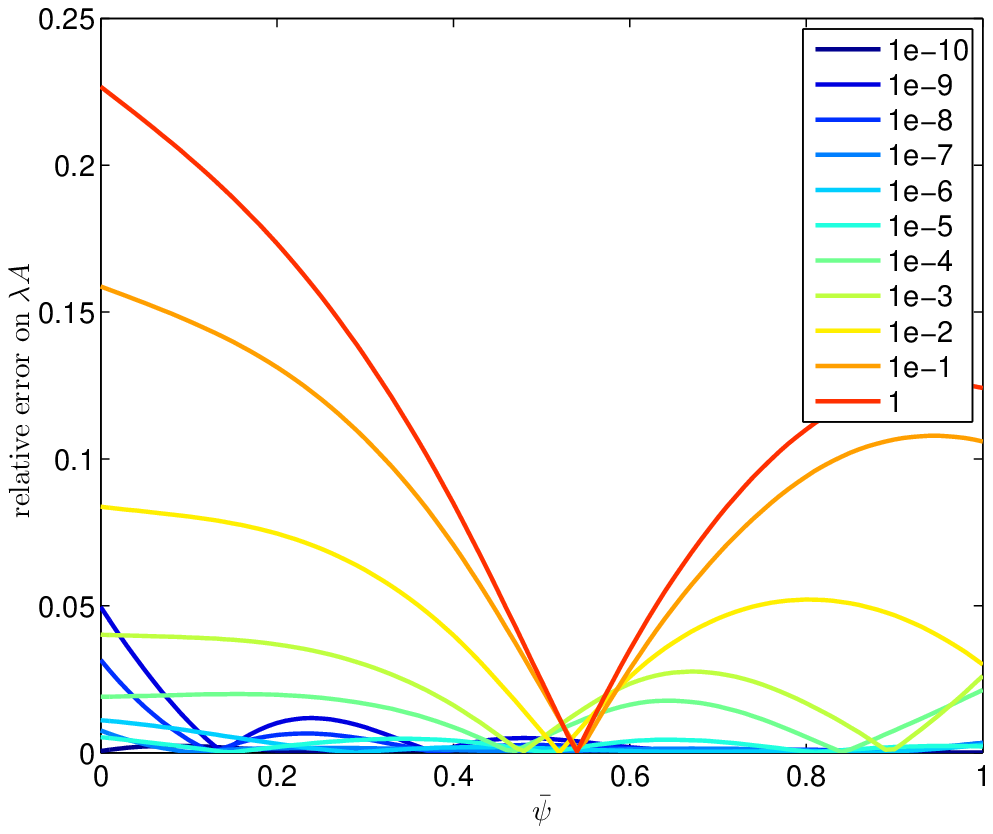} \\  
  \includegraphics[width=6cm]{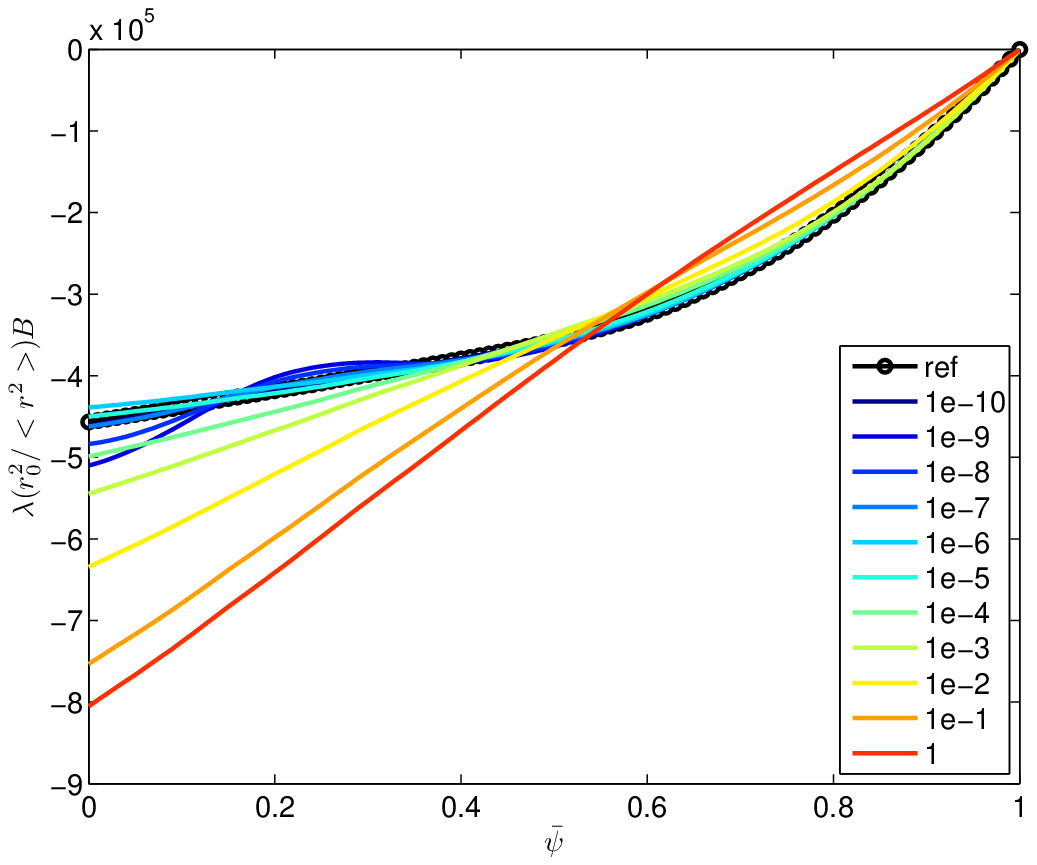} &  
  \includegraphics[width=6cm]{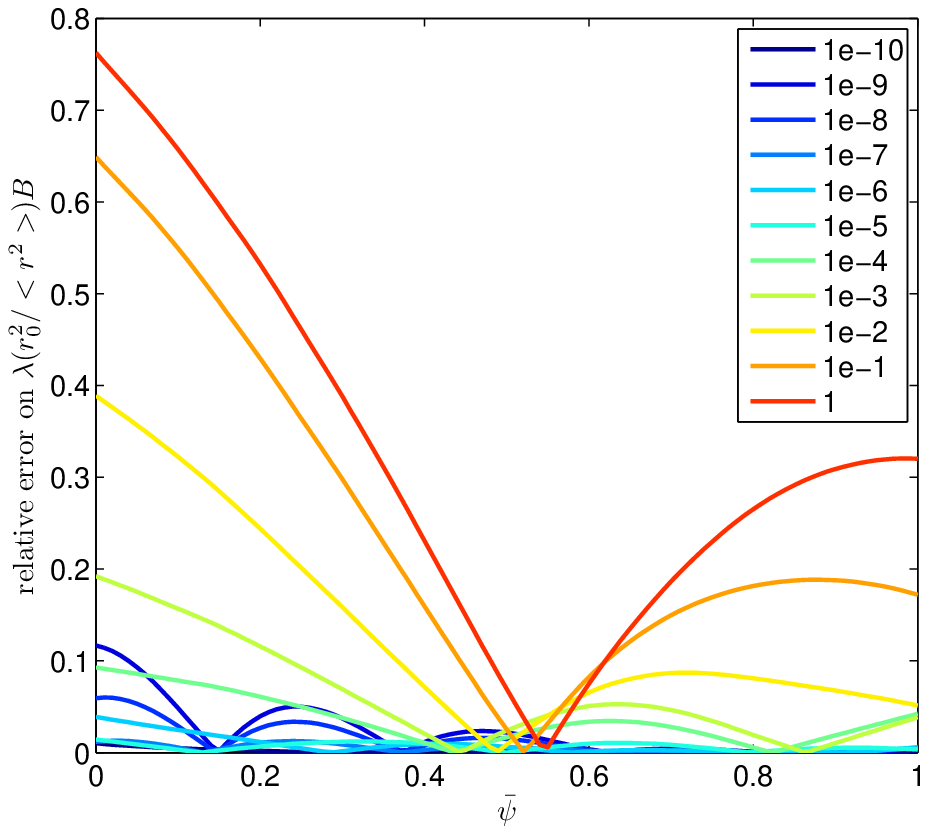}
 \end{tabular}
 \end{center}
 
 \caption{Twin experiment with noise free measurements 
and different regularization parameters $\eps$ ranging from $10^{-10}$ to $1$. 
Left column: identified functions $\lambda A(\p)$ and 
$\lambda R_0^2 \ds <\frac{1}{r^2}> B(\p)$ 
for each different $\eps$ value, 
and the known reference functions. Right column: relative errors. \label{fig:twinAB} }
 \end{figure}

\begin{figure}[!h]
\begin{center}
\begin{tabular}{ll}
  \includegraphics[width=6cm]{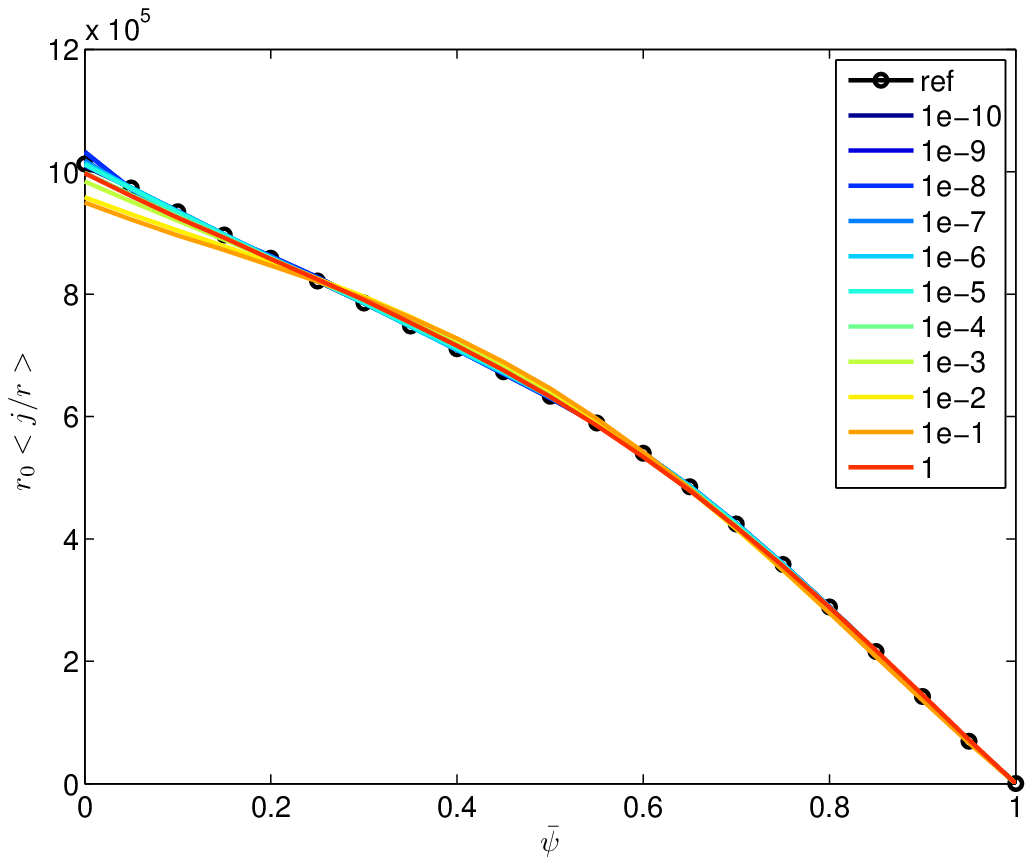} &
  \includegraphics[width=6cm]{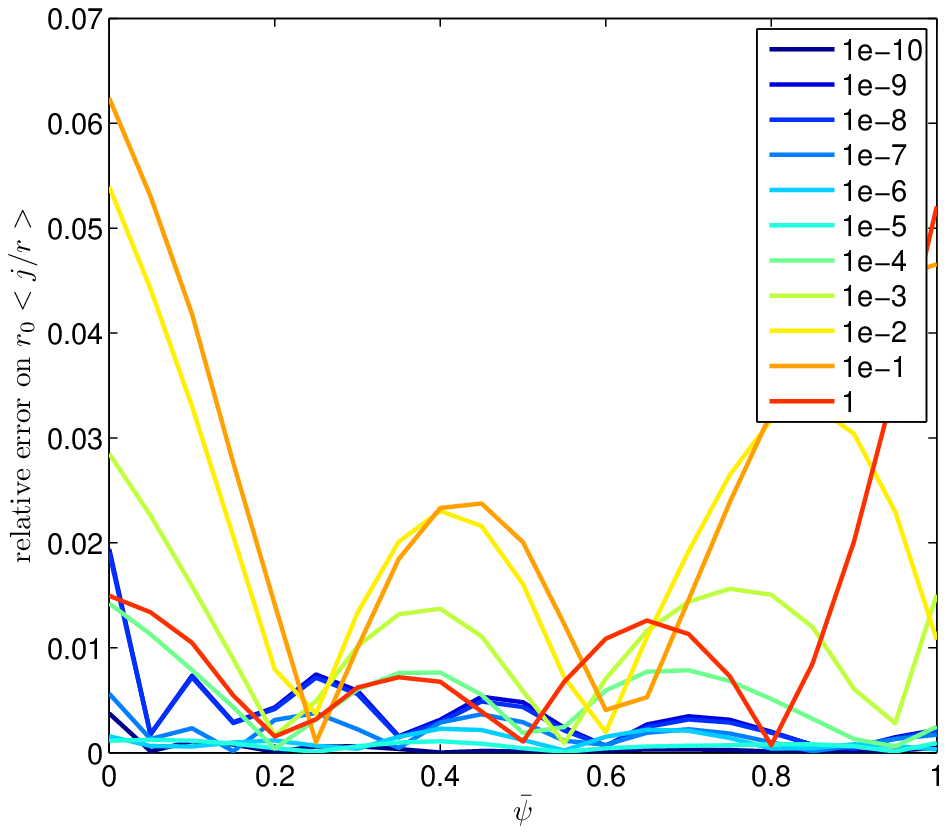}\\ 
  \includegraphics[width=6cm]{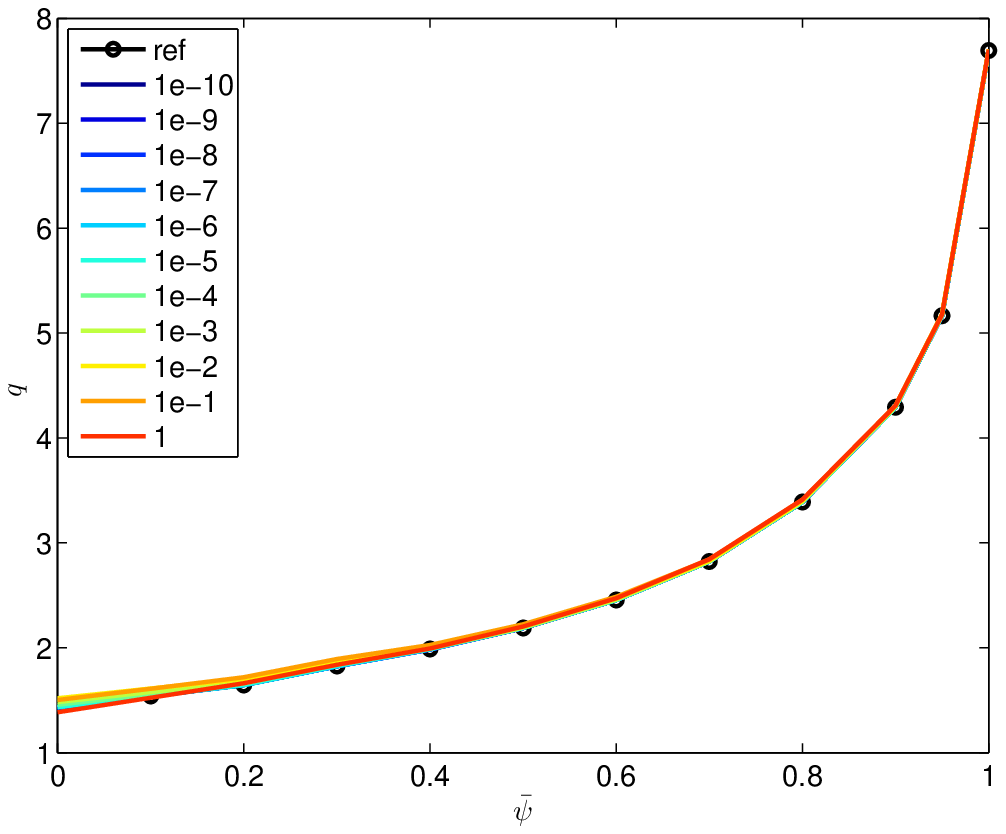} &
  \includegraphics[width=6cm]{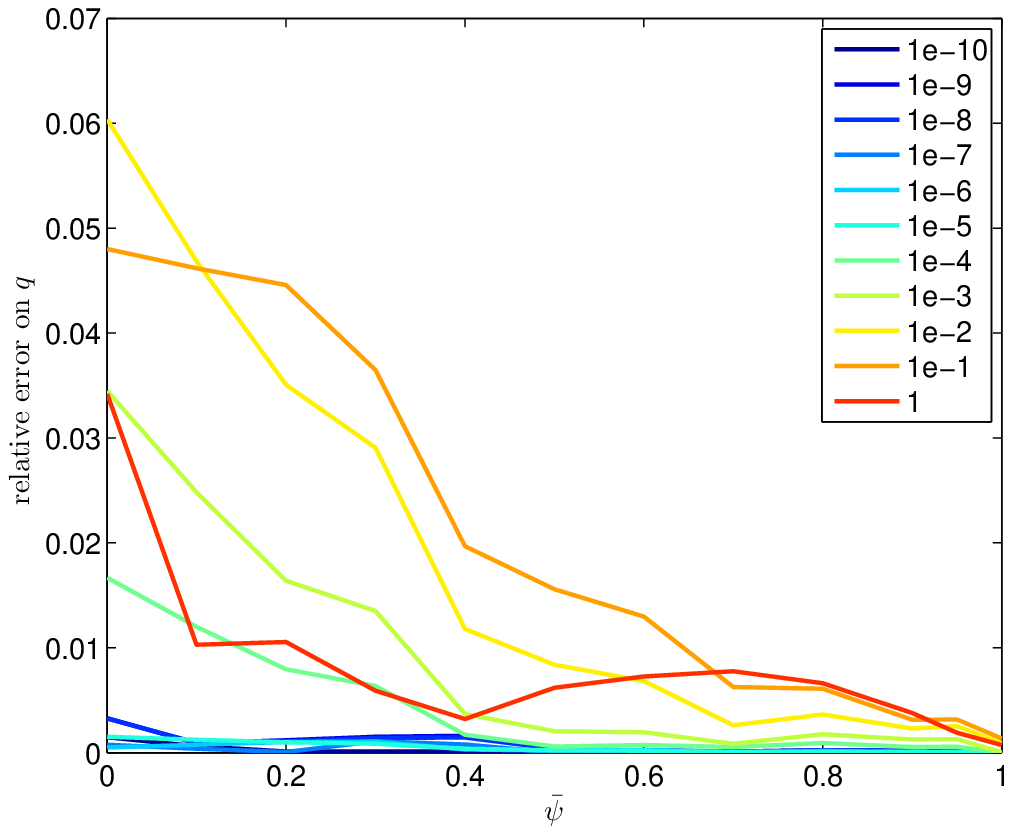} 
 \end{tabular}
 \end{center}
 
 \caption{Twin experiment with noise free measurements 
and different regularization parameters $\eps$ ranging from $10^{-10}$ to $1$. 
Left column: resulting identified averaged current density $R_0<\ds \frac{j(r,\p)}{r}>$, 
safety factor $q$ for each $\eps$ value 
and the corresponding known reference values. Right column: relative errors.
\label{fig:twinjq}}
 \end{figure}

Figure \ref{fig:twinjq} shows an important point. Almost whatever the chosen value of $\eps$ is, 
i.e. whatever the quality of the identification of $A$ and $B$ is, 
the identified averaged current density 
$R_0<\ds \frac{j(r,\p)}{r}>$ as well as the safety factor $q$ are always 
well recovered and the relative errors are one order of magnitude smaller than for 
functions $A$ and $B$. The same kind of observation was made in \cite{Blum:1997} where the identified functions 
$A$ and $B$ seemed to be rather sensitive to perturbations whereas the averaged current density was very stable. 

In Table \ref{tab:residu}, the evolution of the relative
residu on $\psi$, $A$, $B$ and $\lambda$ versus the number of iterations is given. 
It demonstrates numerically the convergence of the algorithm in this case 
where a value of $10^{-6}$ is used as stop condition. The algorithm needs $10$ iterations to converge. 
It is interesting to notice that even though the first guess 
is not particularly well chosen the relative residu on $\psi$ at the second 
iteration has already fallen to $4 \%$. 
In real applications when simulating a whole pulse 
the first guess for the computation of the equilibrium at $t$ 
is the equilibrium computed at $t-\Delta t$ and $2$ iterations are enough to ensure a good 
convergence of the algorithm.

\begin{table}
\caption{\label{tab:residu} Numerical convergence of the algorithm.}
\begin{center}
\begin{tabular}{lllll}
\hline
Iteration $n$ & 
  $\displaystyle{\frac{\|\psi^{n+1}-\psi^n\|}{\|\psi^n\|}}$ &
  $\displaystyle{\frac{\|A^{n+1}-A^n\|}{\|A^n\|}}$ &
  $\displaystyle{\frac{\|B^{n+1}-B^n\|}{\|B^n\|}}$ & 
  $\ds \frac{|\lambda^{n+1}-\lambda^n|}{|\lambda^n|}$ \\
\hline
1&  2.64809 &     6.07599 &    5.3509  & 0.100127      \\  
2&  0.0408642 &   1.19473 &    1.42619 & 9.24968	\\
3&  0.0733385 &   1.83005 &    1.47338  &0.563235	\\
4&  0.0404254 &   0.884617 &   1.0359 & 0.108107	\\
5&  0.00539736 &  4.79091 &    4.37571 & 0.826455	\\
6&  0.000349811 & 0.127626 &   0.180449 &0.0889022	\\
7&  1.58606e-05 & 0.0262942 &  0.0246657 & 0.0263	\\
8&  5.67036e-06 & 0.00294791 & 0.0024952 & 0.00315952	\\
9&  1.4533e-06 &  0.000339986 &0.000273055& 0.000362224\\
10& 6.19066e-07 & 6.41923e-05 &6.51076e-05& 6.29838e-05\\
\hline
\end{tabular}
\end{center}
\end{table}

\subsection{Twin experiment with noisy magnetic measurements}

\begin{figure}[!h]
\begin{center}
\begin{tabular}{ll}
  \includegraphics[width=6cm]{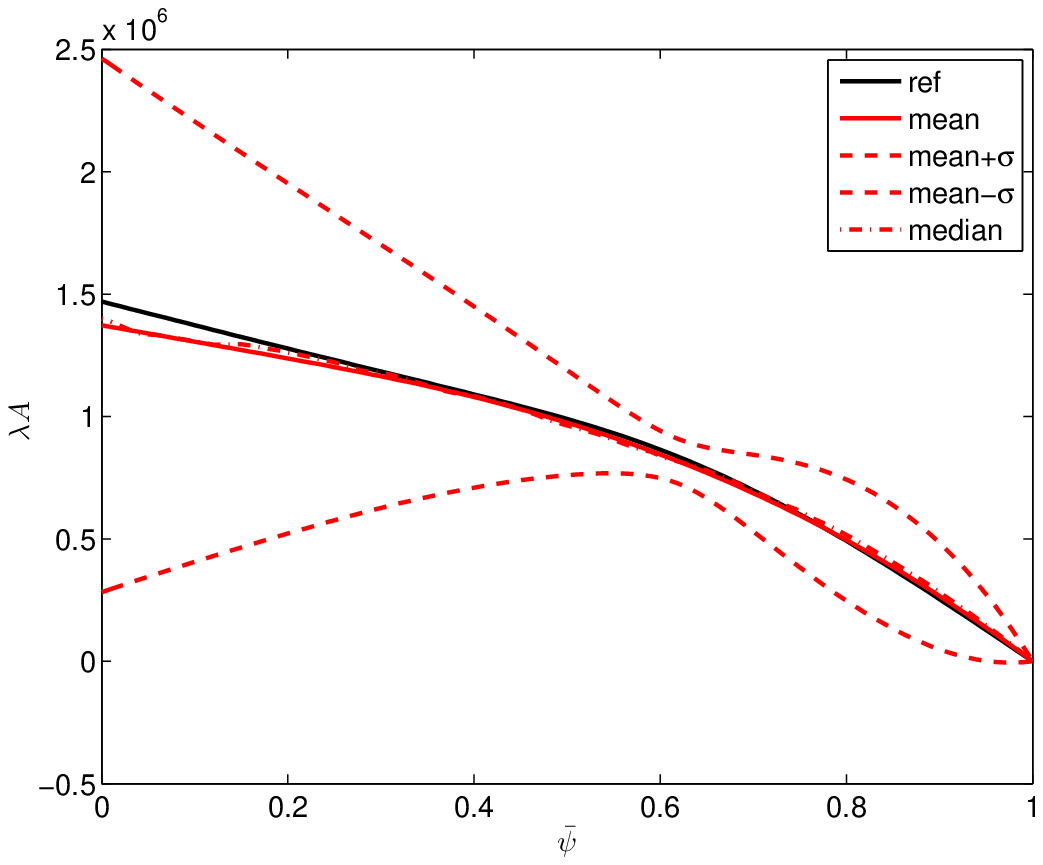} & \includegraphics[width=6cm]{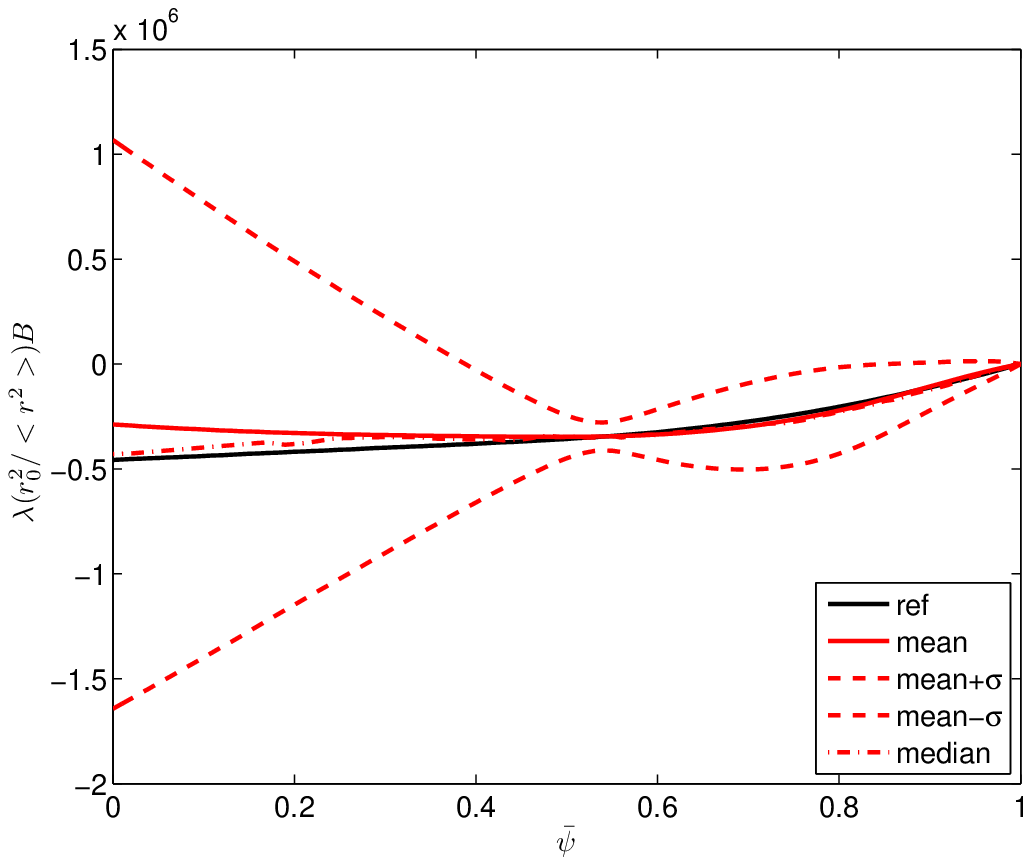}\\  
  \includegraphics[width=6cm]{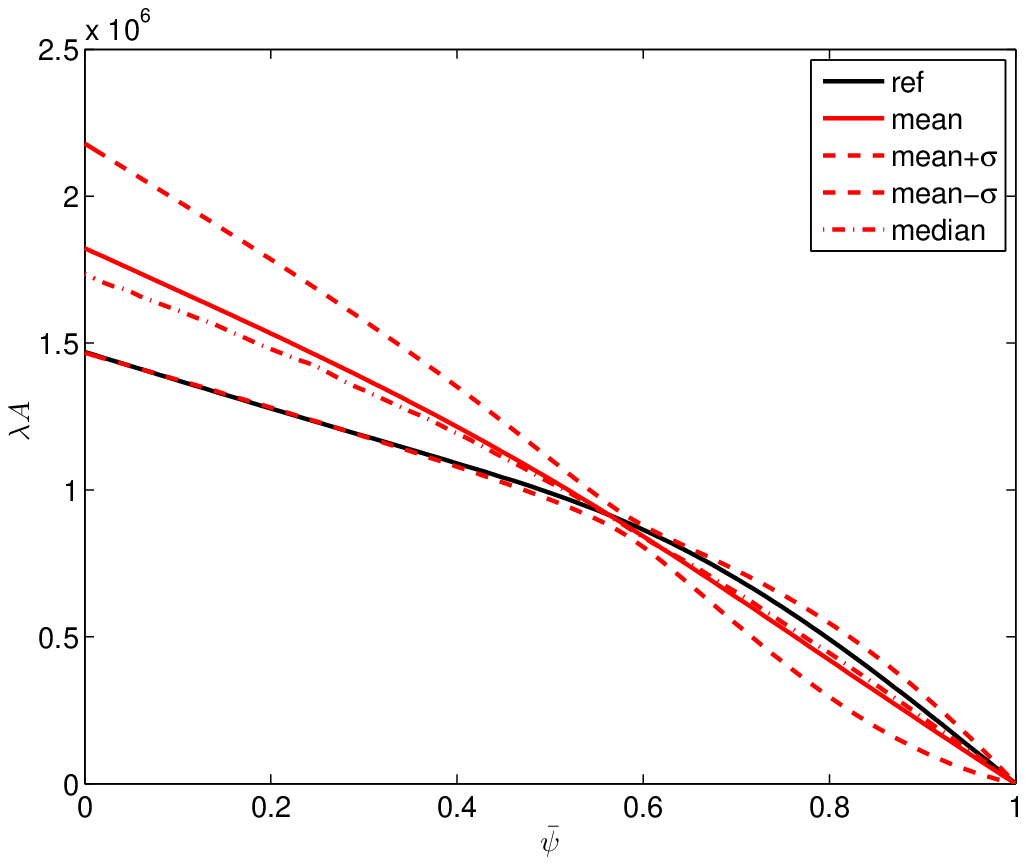} &  \includegraphics[width=6cm]{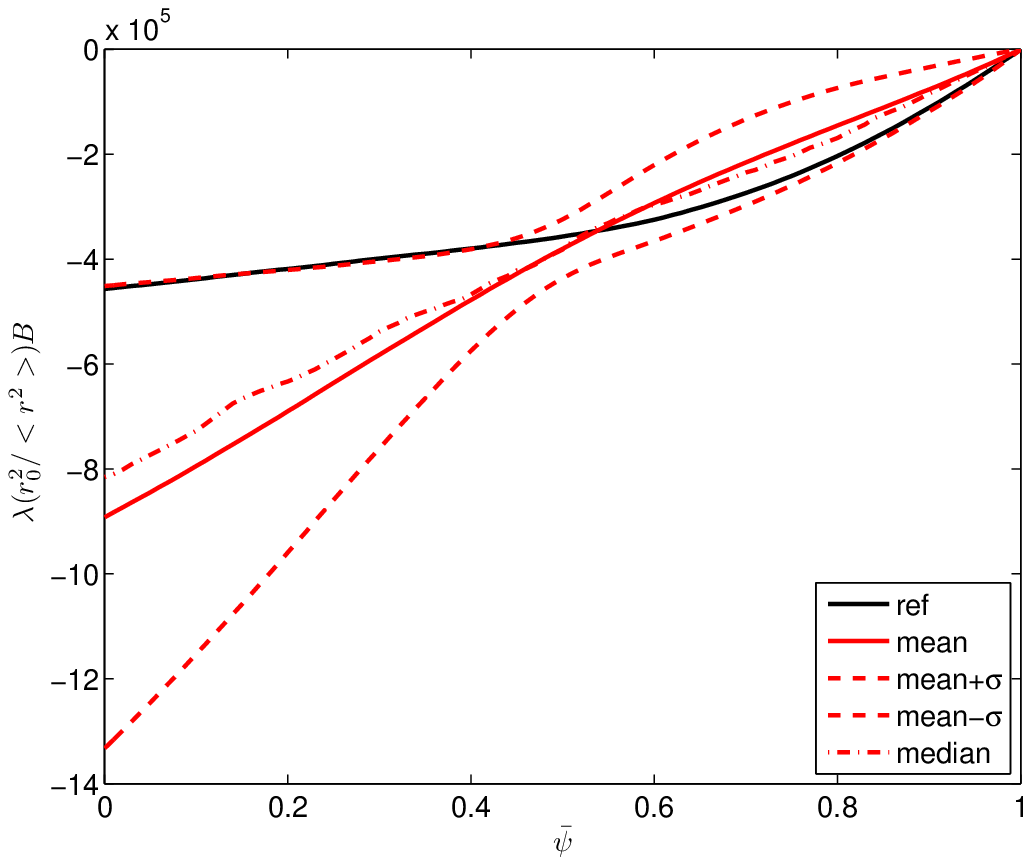} \\
  \includegraphics[width=6cm]{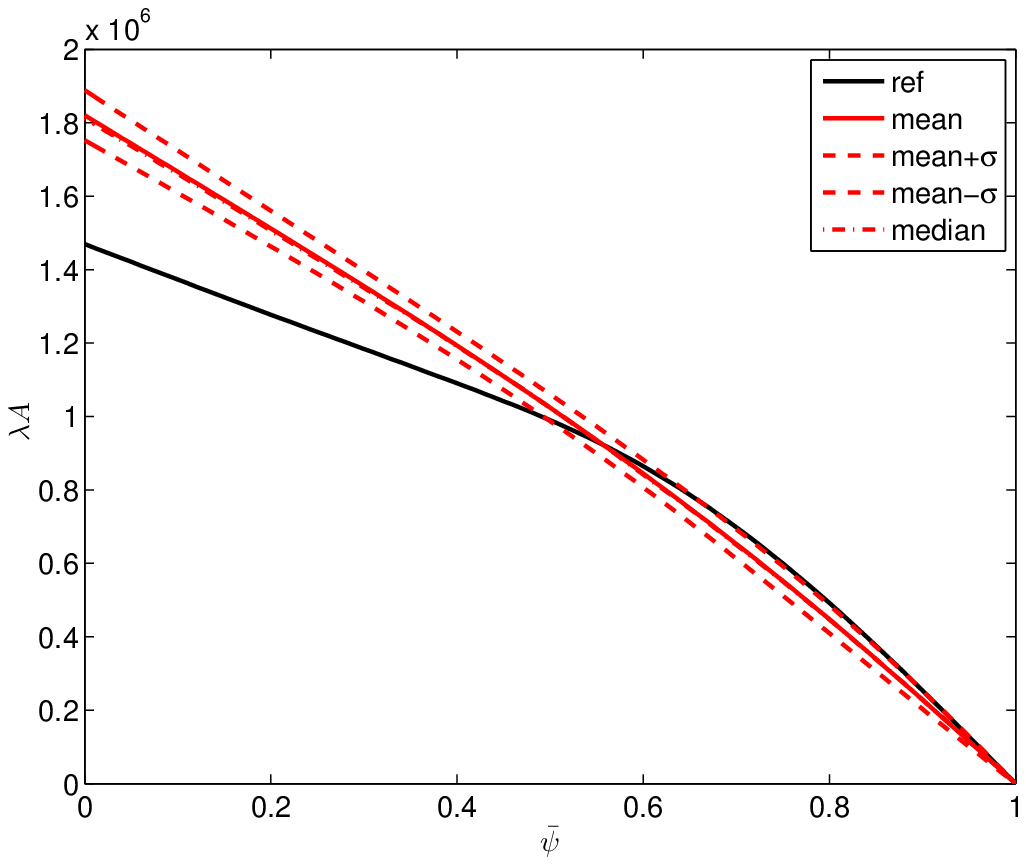} &  \includegraphics[width=6cm]{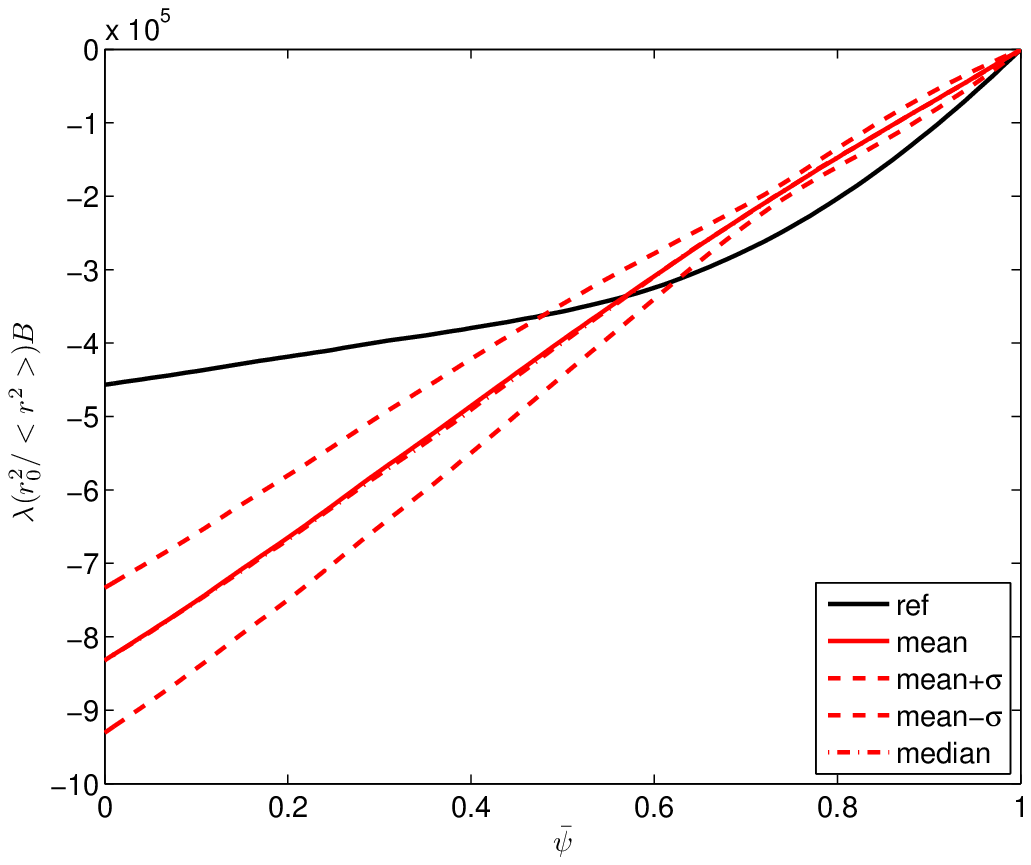} 
\end{tabular}
\end{center}
\caption{Statistical results of the identification experiments with noisy magnetic measurements. 
Row 1: $\eps=10^{-2}$, row 2: $\eps=10^{-1}$, row 3 $\eps=1$. Column 1: function $\lambda A(\p)$ 
and column 2: $\lambda R_0^2 \ds <\frac{1}{r^2}> B(\p)$. For each function the reference value from which 
the unperturbed measurements were computed is given in black and the mean identified function in red. 
The mean $\pm$ standard deviation functions are shown in dashed red. \label{fig:twinmeannoise}}
\end{figure}

\begin{figure}[!h]
\begin{center}
\begin{tabular}{lll}
  \includegraphics[width=6cm]{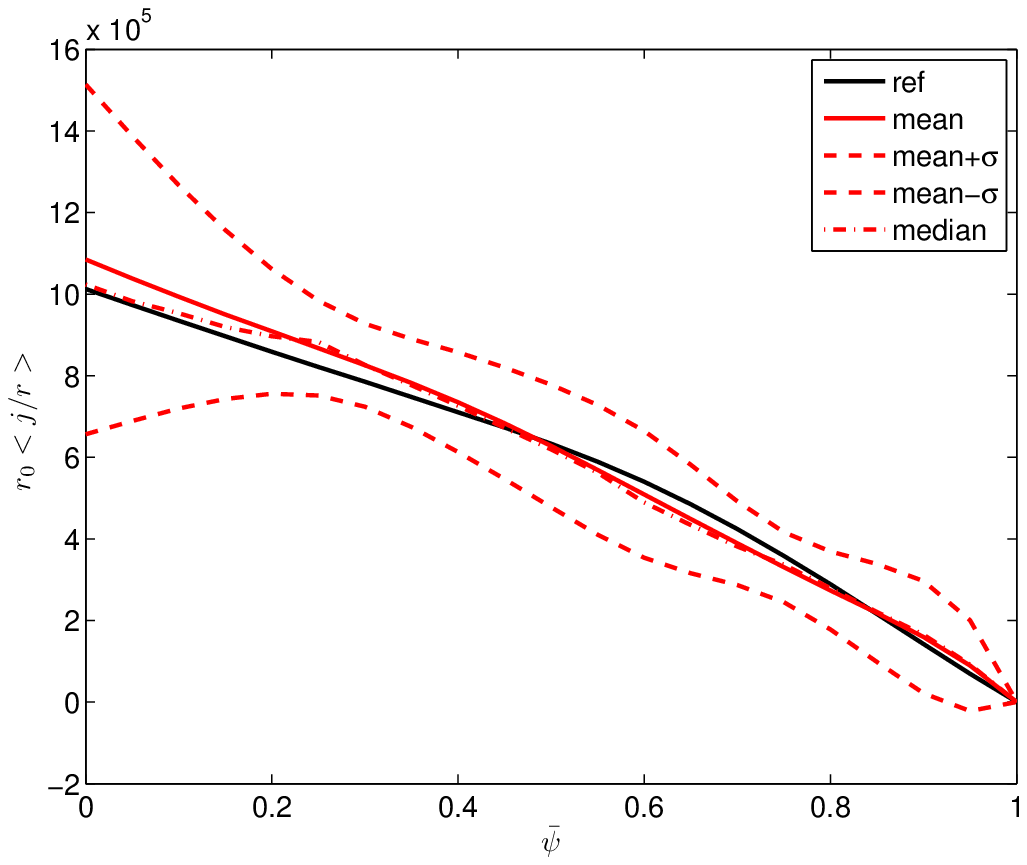} & \includegraphics[width=6cm]{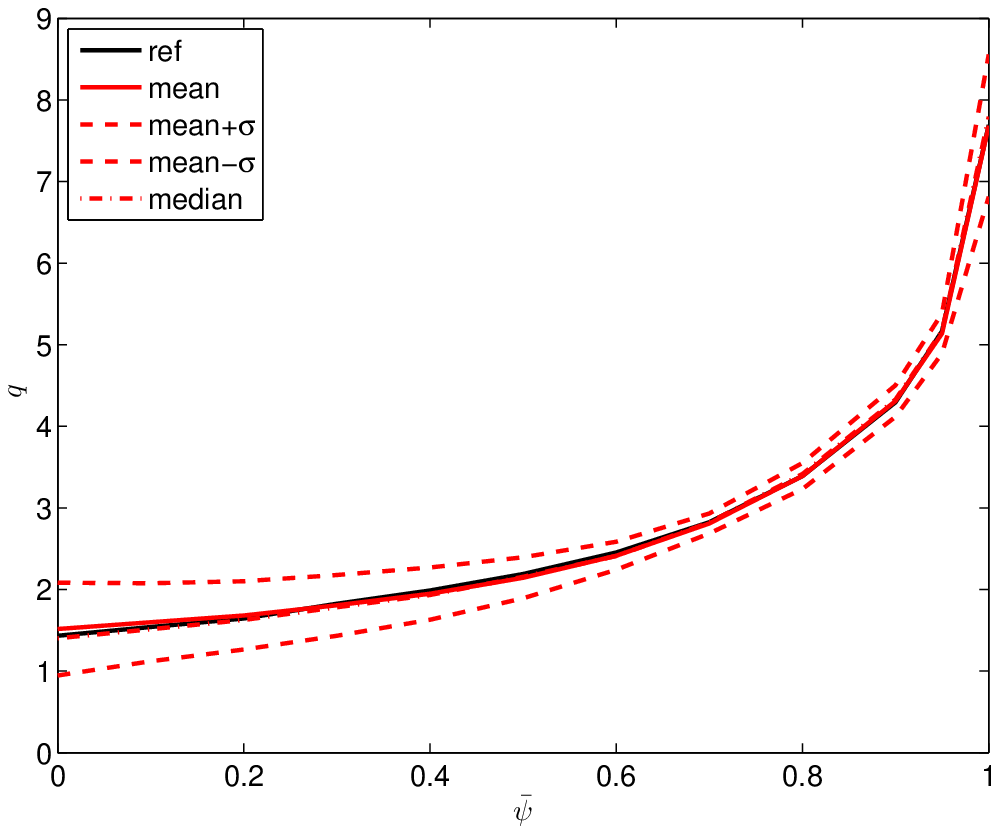} \\
  \includegraphics[width=6cm]{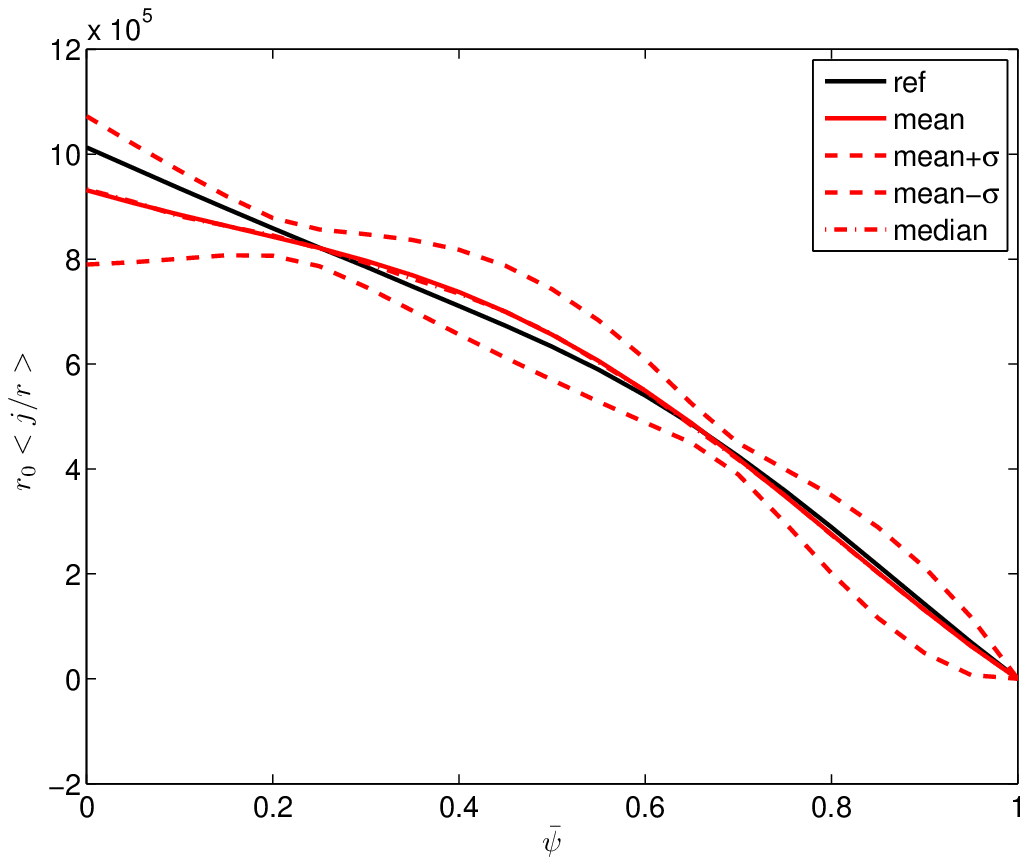} &  \includegraphics[width=6cm]{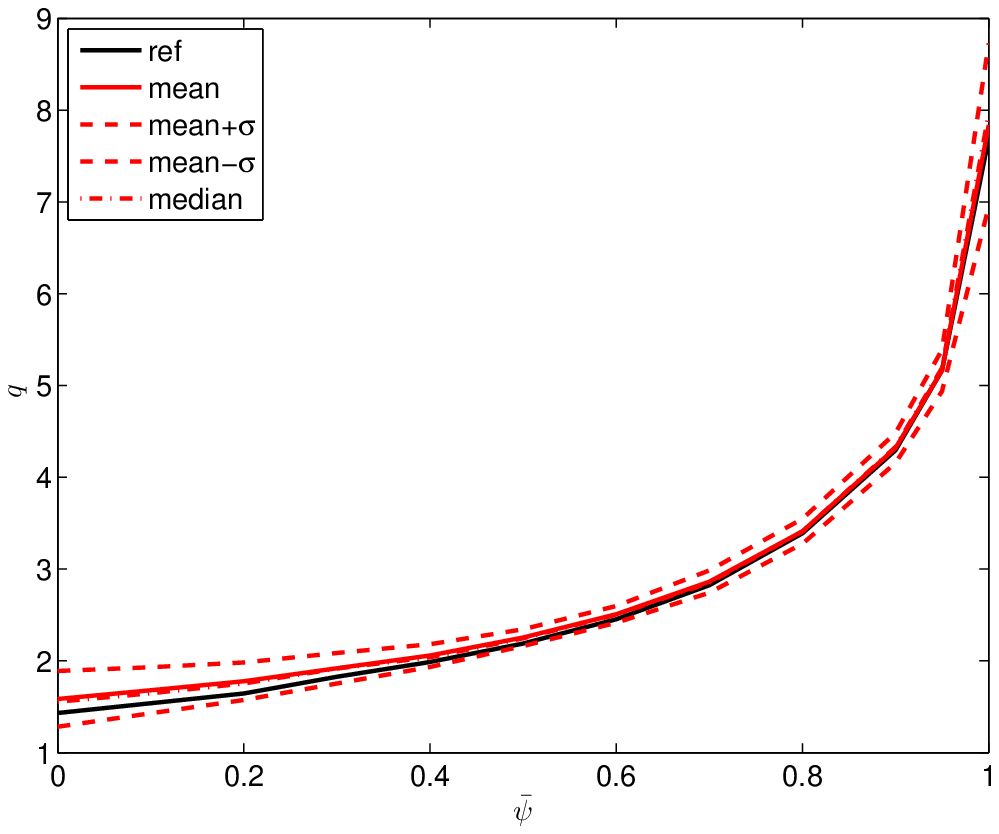} \\  
  \includegraphics[width=6cm]{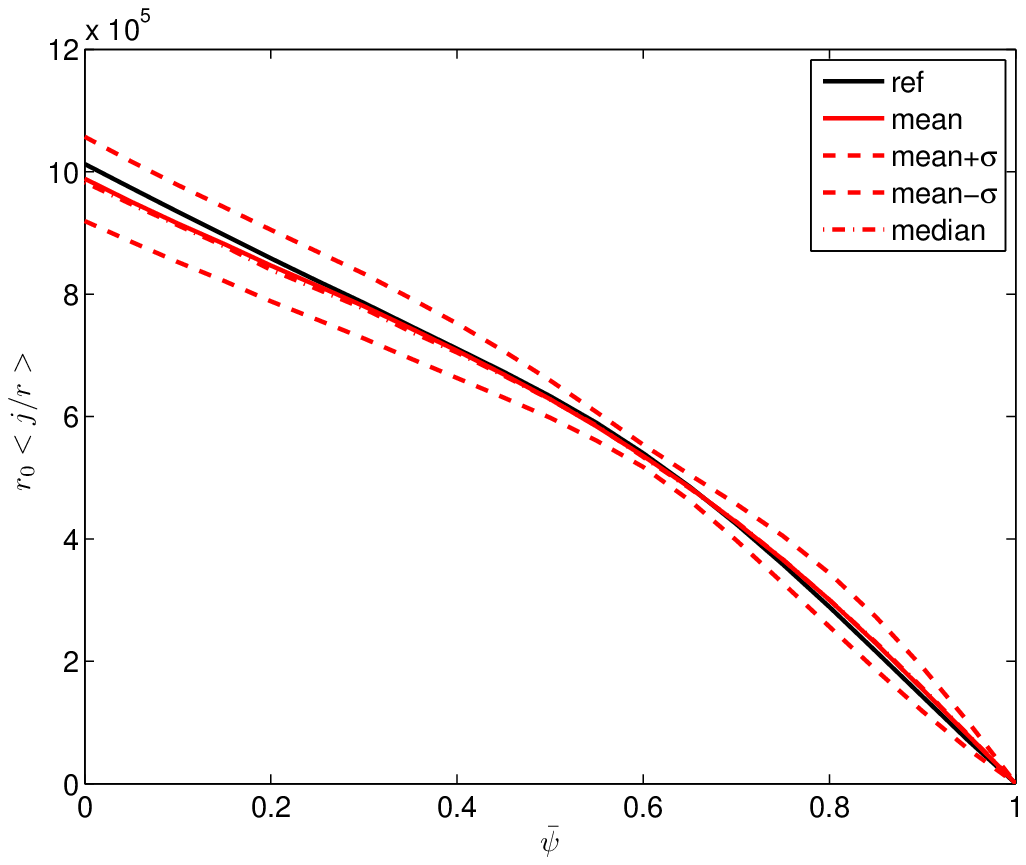} & \includegraphics[width=6cm]{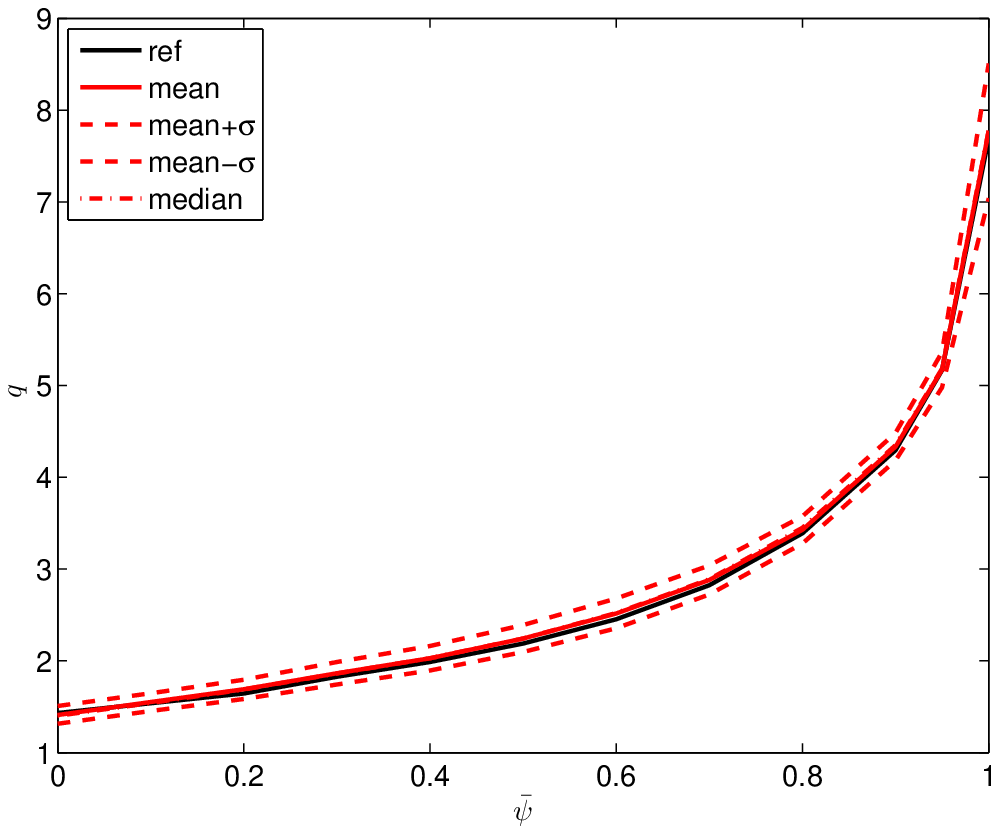} 
\end{tabular}
\end{center}
\caption{Statistical results of the identification experiments with noisy magnetic measurements. 
Row 1: $\eps=10^{-2}$, row 2: $\eps=10^{-1}$, row 3 $\eps=1$. Column 1: $ R_0<\ds \frac{j(r,\p)}{r}>$, 
and column 2: safety factor $q$.
For each function the reference value is given in black and the mean identified function in red. 
The mean $\pm$ standard deviation functions are shown in dashed red. \label{fig:twinmeannoise2}}
\end{figure}

Figures \ref{fig:twinmeannoise} and \ref{fig:twinmeannoise2} 
show the results of the same type of numerical experiment but with noisy measurements. 
Each magnetic input, $m$ representing either $\psi$ 
or $\ds \frac{1}{r}\frac{\partial \psi}{\partial n}$
at a point of the domain boundary $\Gamma$, is perturbed 
with a one percent noise normally distributed, 
$m_{\eta}=m+\eta$ with $\eta \sim N(m,0.01 m)$. 
For each chosen value of the regularization parameter the algorithm is run 200 times 
with measurements randomly perturbed as above. Then for each function $\lambda A$, 
$\lambda \ds R_0^2 \ds <\frac{1}{r^2}>B$, $R_0 <\ds \frac{j(r,\p)}{r}>$ and $q$, a mean function 
and a standard deviation function are computed.

In comparison with the noise free case the regularization parameter 
needs to be significantly increased to values of at least $\eps=10^{-2}$ 
and for a safer convergence of the algorithm to $\eps=10^{-1}$. 
For smaller values the algorithm either does not converge 
or gives very oscillating identified functions.

The mean error on the reconstructed functions is always smaller 
in the interval $\p \in [0.5,1]$ than in 
the interval $[0,0.5]$. This is due to the fact that magnetic measurements are external to the plasma and 
do not provide enough information to properly reconstruct the functions in the innermost part of the plasma.

As $\eps$ increases the variability or the standard deviation on the identified functions decreases. 
With small $\eps$ the algorithm can find very different functions 
depending on the perturbations of 
the measurements. With $\eps=10^{-2}$ the variability in the identified 
functions $A$ and $B$ 
is strong however the mean identified functions are close to the exact reference ones. On the other hand 
with $\eps=1$ the variability of the identified functions 
is strongly reduced but they are quite 
different from the exact reference functions in the interval $[0,0.5]$.

It is worth noticing that in all cases the resulting 
safety factor $q$ and averaged current density 
$R_0<\ds \frac{j(r,\p)}{r}>$ are well recovered. 
The remark of the preceding section on the identifiability of the averaged current 
density still holds: 
it is quite well recovered even if functions $A$ and $B$ 
taken separately are not well identified. 
The mean error on the current density profile is almost always smaller than the mean errors on functions 
$A$ and $B$. Moreover this error does not change very much between the different cases and particularly 
between the $\eps=10^{-1}$ and the $\eps=1$ cases. 
This implies that for a large interval of $\eps$ the value of the part of the cost function 
related to magnetic measurements $J_0$ is almost constant. 
Therefore it is difficult to find an optimal value for the regularization parameter. 
For example the L-curve method \cite{Hansen:1999} for the determination of 
the regularization parameter can hardly be used and 
gives some results which are not very reliable since the 
L-curves are not well behaved and the location of the corner is not clear. 
The "L" is an almost vertical line. This is due to the fact that, in a large interval of $\eps$ 
values, an increase in $\eps$ implies a important decrease 
in the regularization term $\frac{1}{2} (u^*(\eps))^T \Lambda u^*(\eps)$ 
but does not lead to a significative increase in the misfit term $J_0(u^*(\eps))$.

\clearpage
\newpage

\subsection{Twin experiment with noisy magnetic, interferometric 
and polarimetric measurements}

In this last twin experiment, interferometric and polarimetric 
measurements are also used. At first a reference density profile, $n_e(x)$ is prescribed 
point by point on $[0,1]$,
as well as the same reference $A$ and $B$ functions as in the previous twin experiments. 
Then similar to the preceding section 
the equilibrium is computed from given 
Dirichlet boundary condition. 
A set of artificial magnetic, interferometric 
and polarimetric measurements is generated. 
Finally several twin experiments with a $1\%$ noise are performed and some statistics are computed. 
The weights related to interferometric and polarimetric measurements in the cost function are 
defined as 
\begin{itemize}
\item $w_k^{polar}=\ds \frac{1}{\sqrt{N_c} \sigma^{polar}}$, with $\sigma^{polar}=10 ^{-1}$ radians
\item $w_k^{inter}=\ds \frac{1}{\sqrt{N_c} \sigma^{inter}}$, with $\sigma^{inter}=10 ^{18}$ $m^{-3}$
\end{itemize}

\begin{figure}[!h]
\begin{center}
\begin{tabular}{lll}
  \includegraphics[width=7cm]{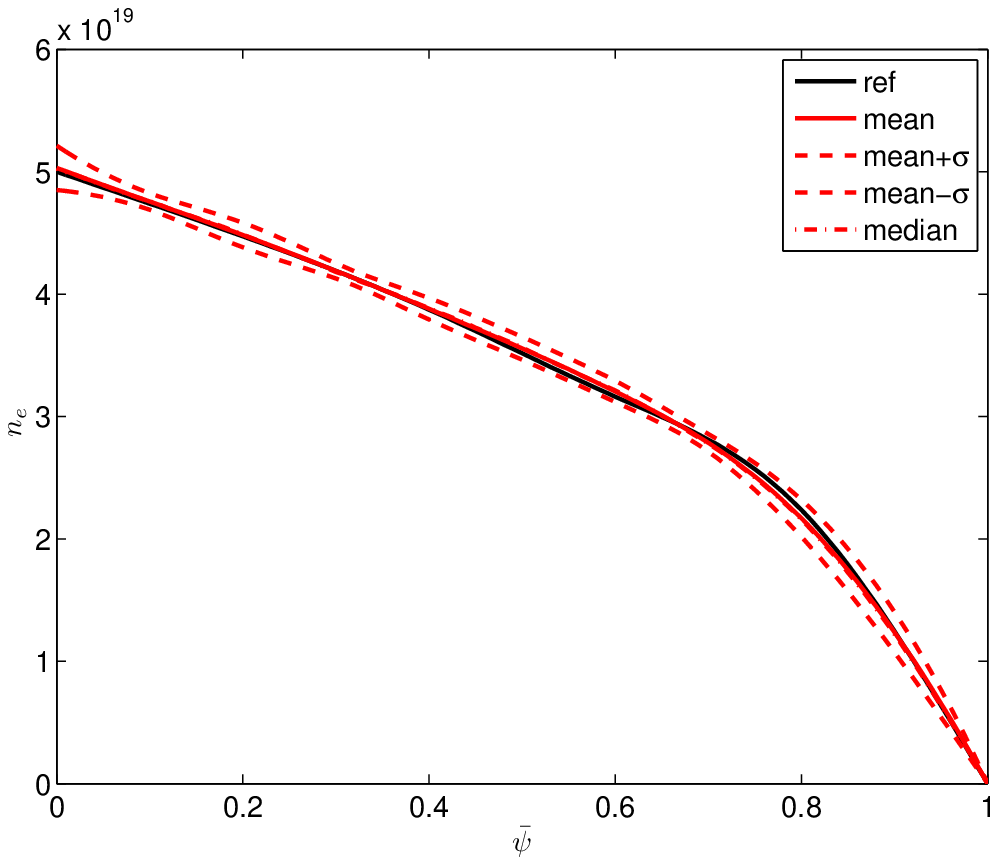} 
 \end{tabular}
 \end{center}
 \caption{Statistical results for the identification of the density function $n_e$ 
with noisy interferometric measurments. \label{fig:twinmeannoisene}}
 \end{figure}

The determination of the regularization parameter for the density function $n_e$ is far less 
a problem than for functions $A$ and $B$ since for example the L-curve method works quite well in this case 
(see Fig. \ref{fig:LcurveNe} in the next Section) and the $n_e$ function is well recovered as shown on 
Fig. \ref{fig:twinmeannoisene}. The regularization parameter for the density function is set to 
$\eps_{ne}=10^{-2}$. 

\begin{figure}[!h]
\begin{center}
\begin{tabular}{lll}
  \includegraphics[width=6cm]{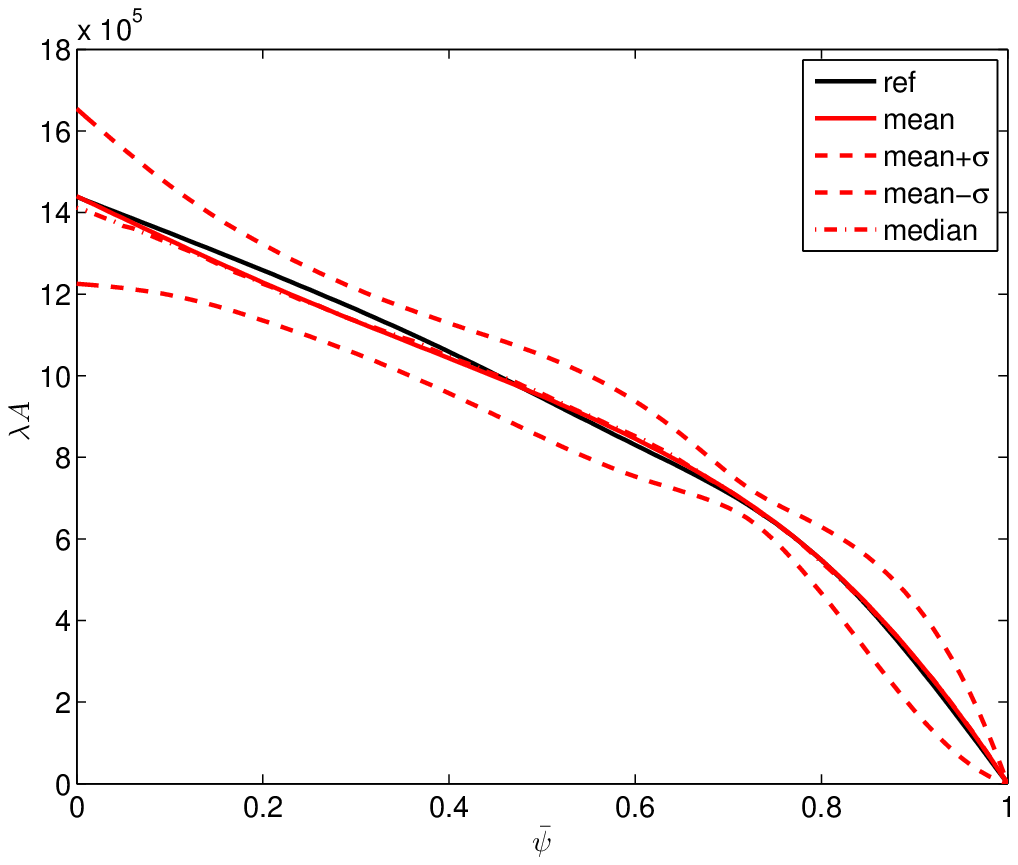} & \includegraphics[width=6cm]{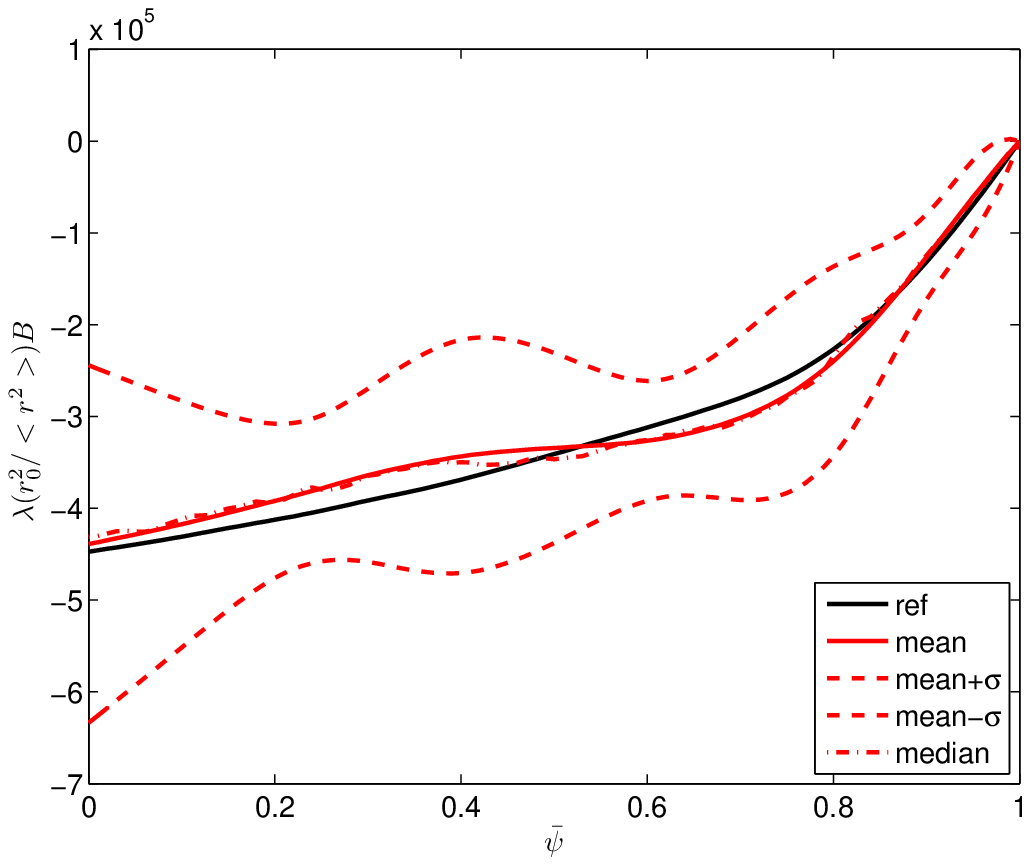}\\
\includegraphics[width=6cm]{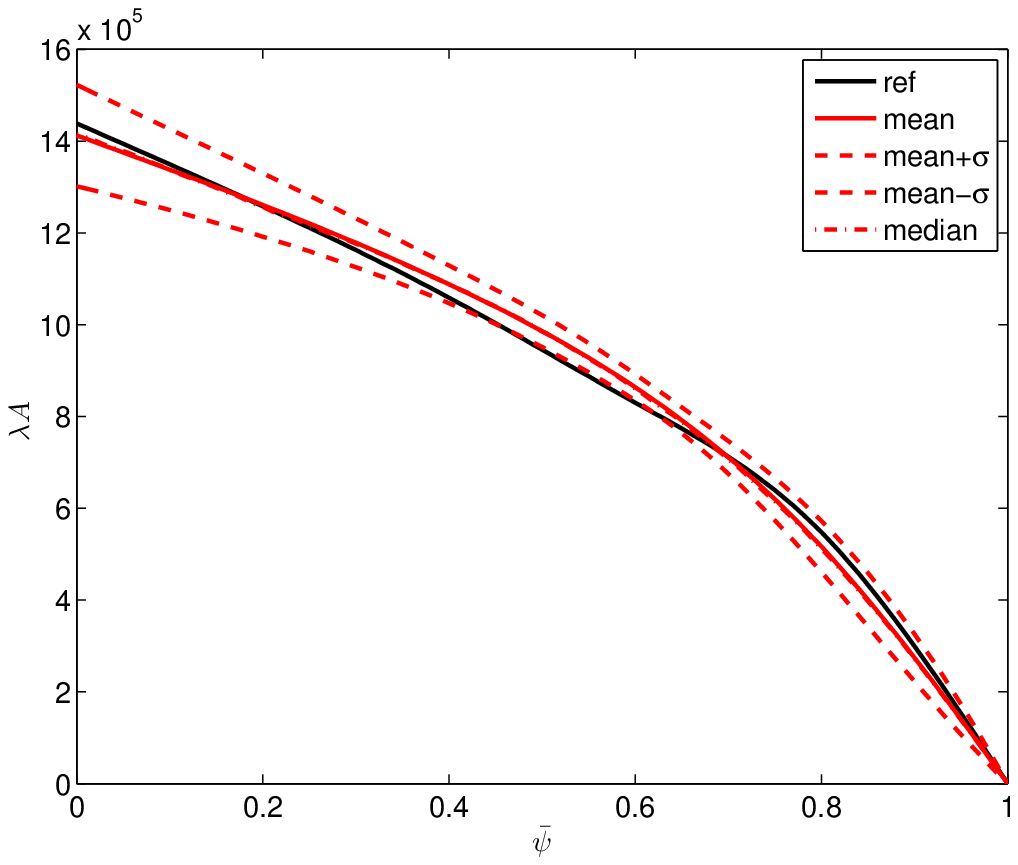} & \includegraphics[width=6cm]{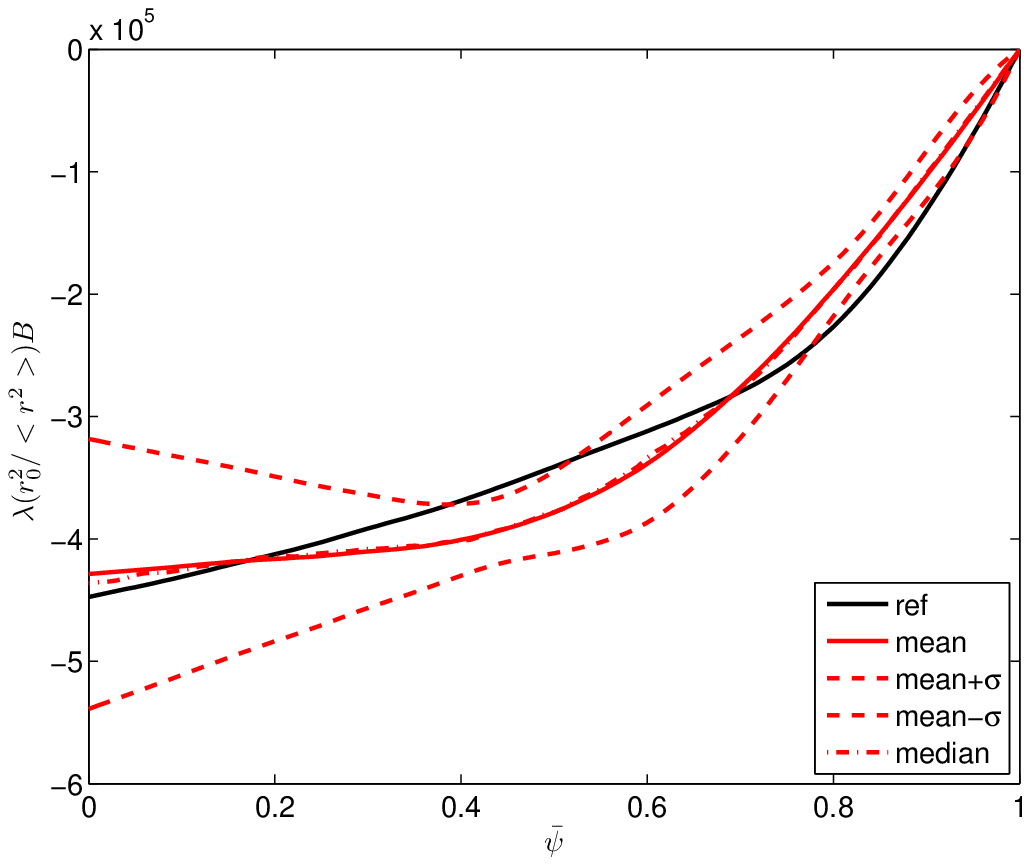}\\
\includegraphics[width=6cm]{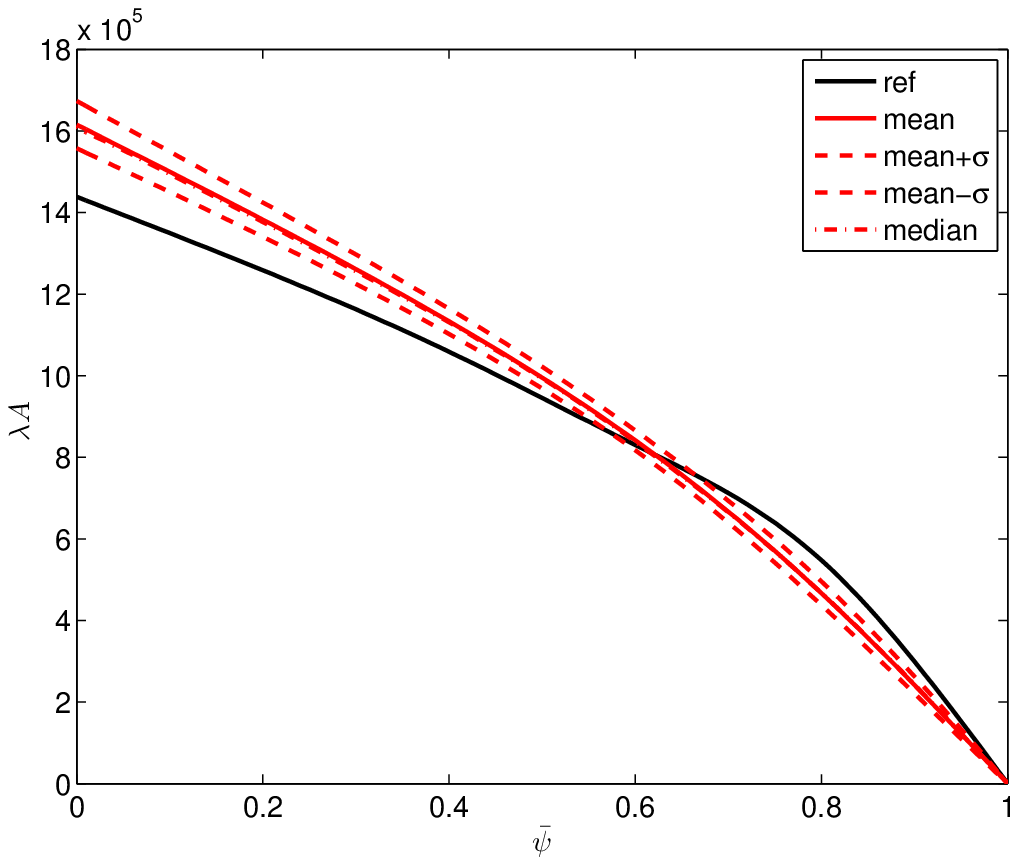} & \includegraphics[width=6cm]{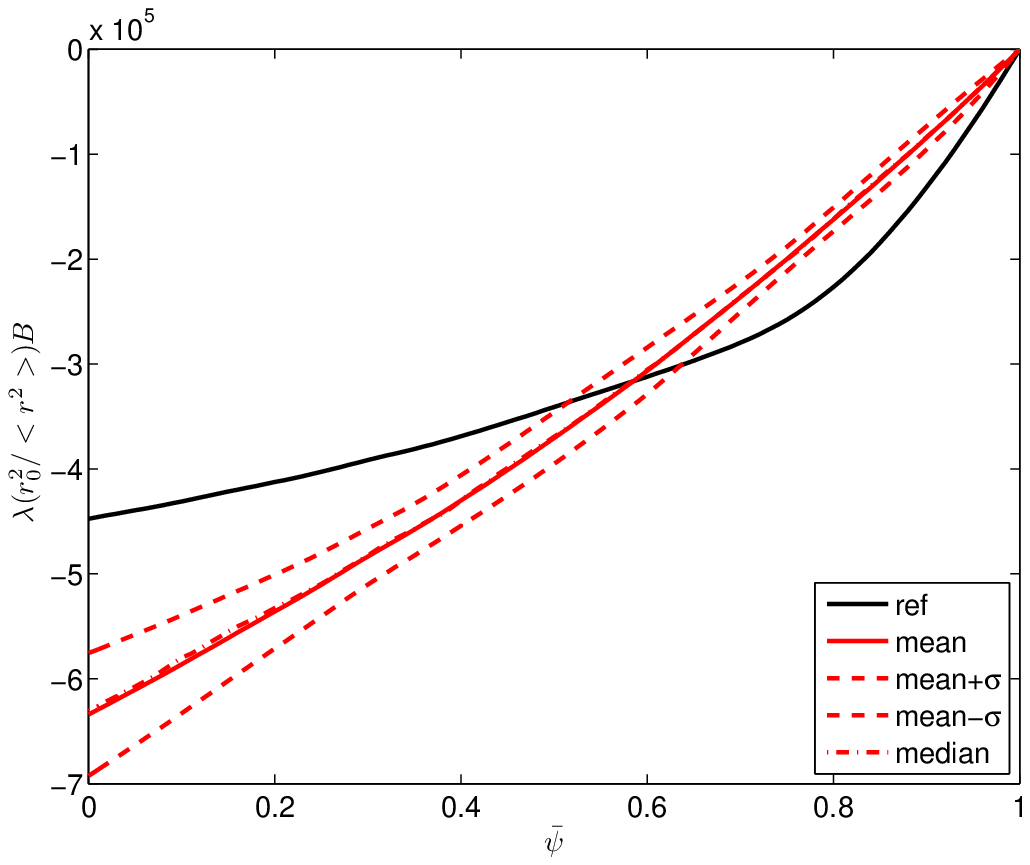} 
 \end{tabular}
 \end{center}
 \caption{
Statistical results of the identification experiments with noisy measurements 
(magnetics, interferometry and polarimetry). 
Row 1: $\eps=10^{-2}$, row 2: $\eps=10^{-1}$, row 3 $\eps=1$. Column 1: function $\lambda A(\p)$, 
and column 2: $\lambda R_0^2 \ds <\frac{1}{r^2}> B(\p)$.
For each function the reference value from which the unperturbed measurements were computed is given 
in black and the mean identified function in red. 
The mean $\pm$ standard deviation functions are shown in dashed red. \label{fig:twinmeannoisepolar}}
 \end{figure}

\begin{figure}[!h]
\begin{center}
\begin{tabular}{lll}
\includegraphics[width=6cm]{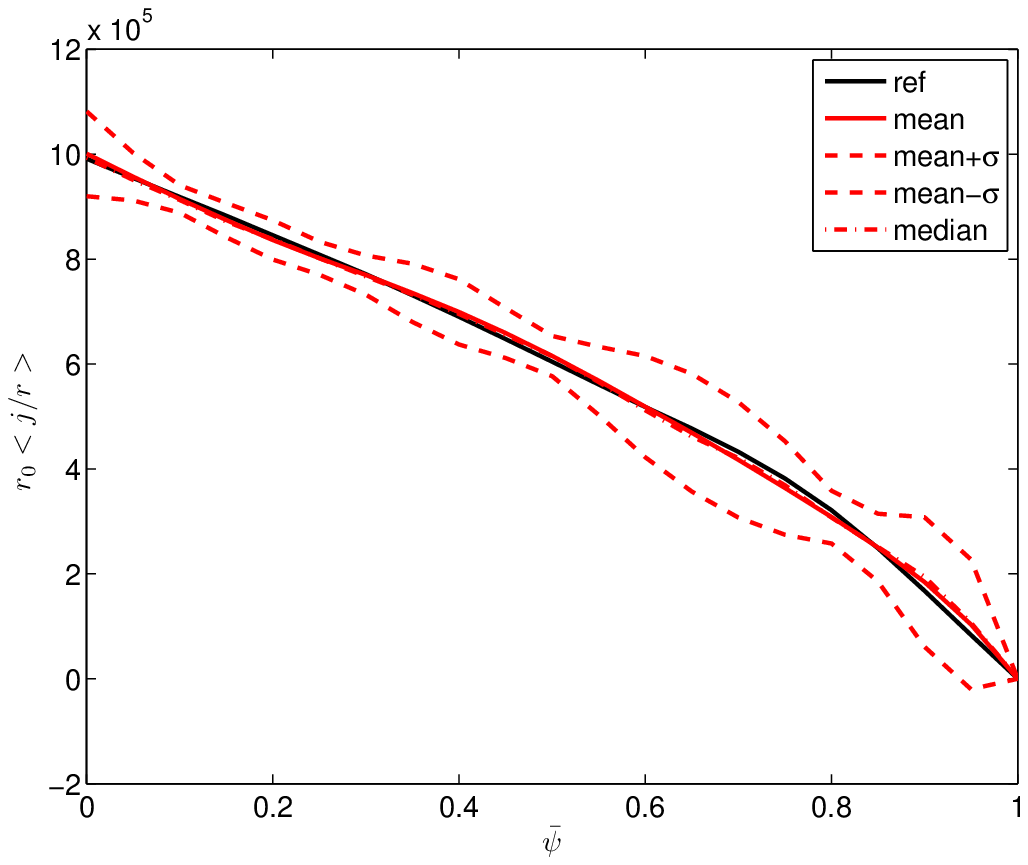}&\includegraphics[width=6cm]{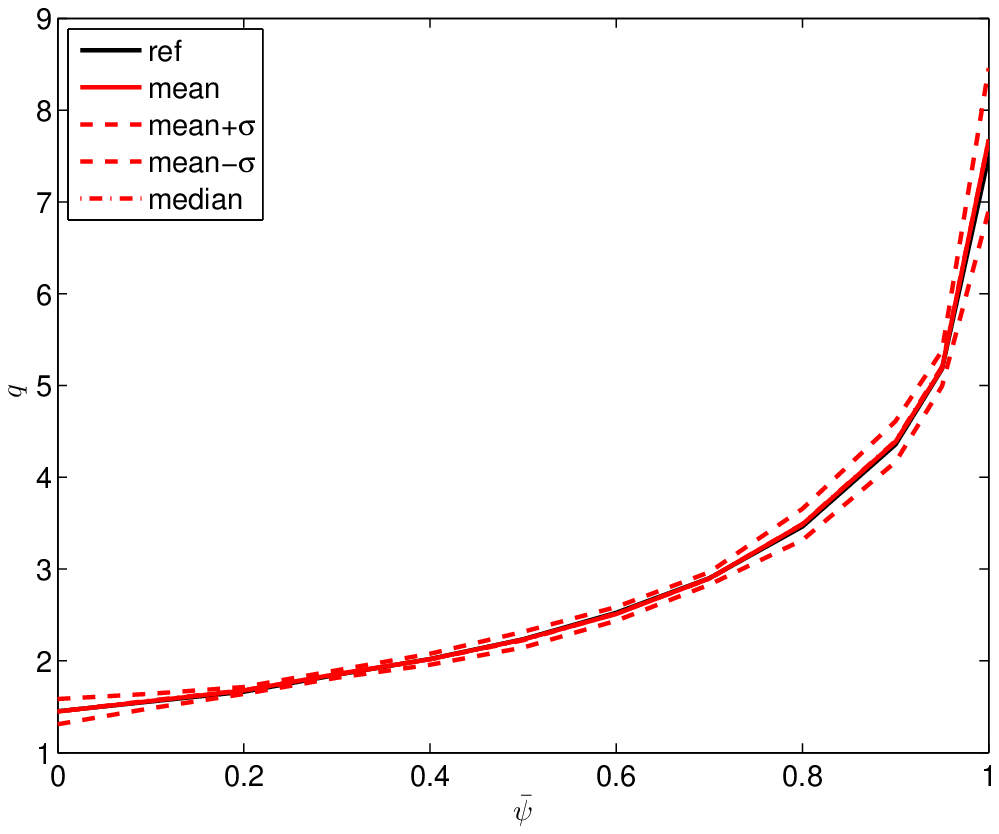}\\
\includegraphics[width=6cm]{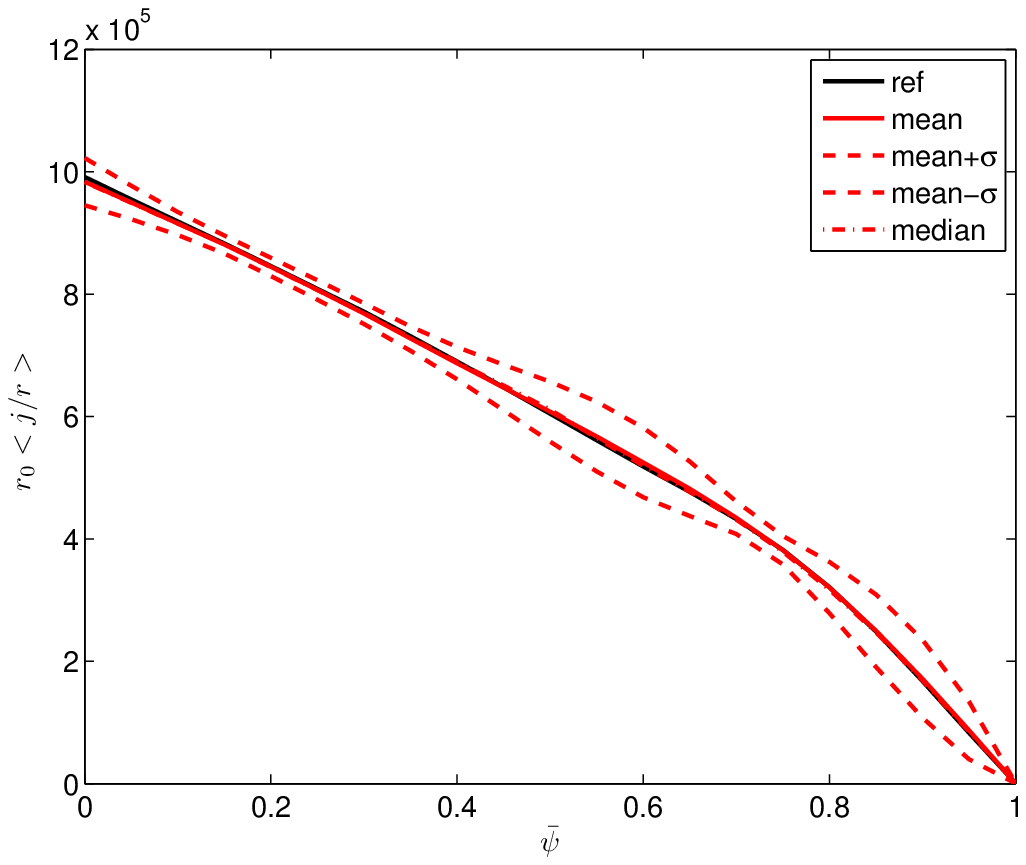}&\includegraphics[width=6cm]{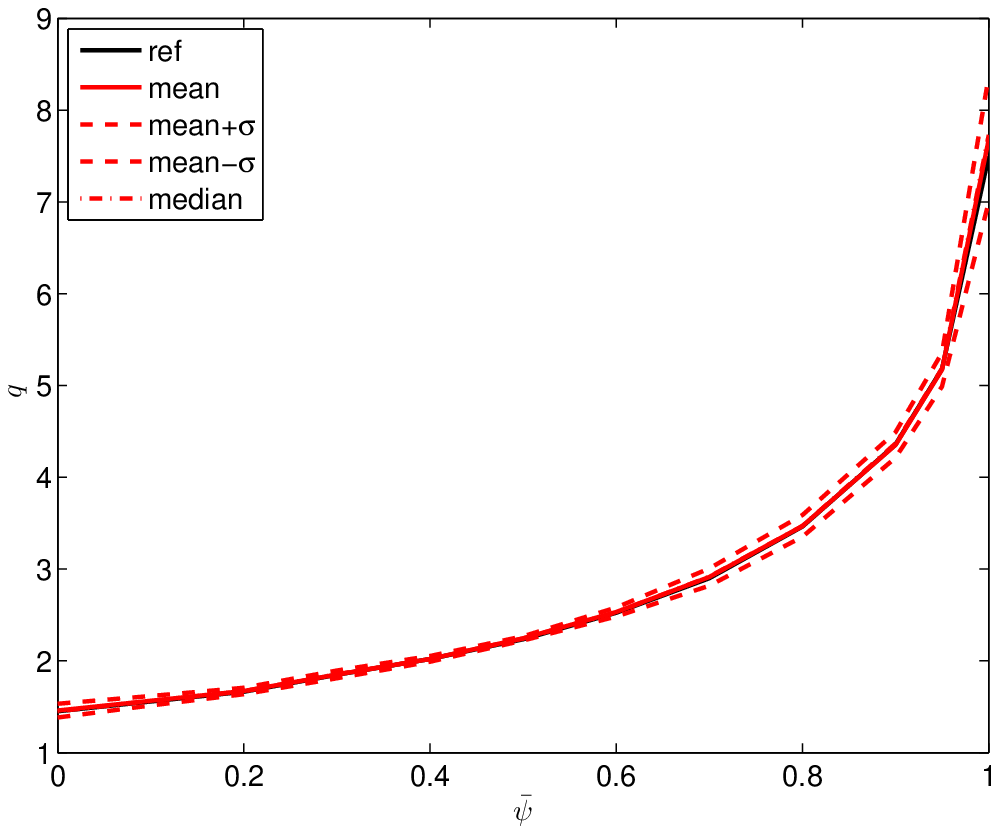}\\  
\includegraphics[width=6cm]{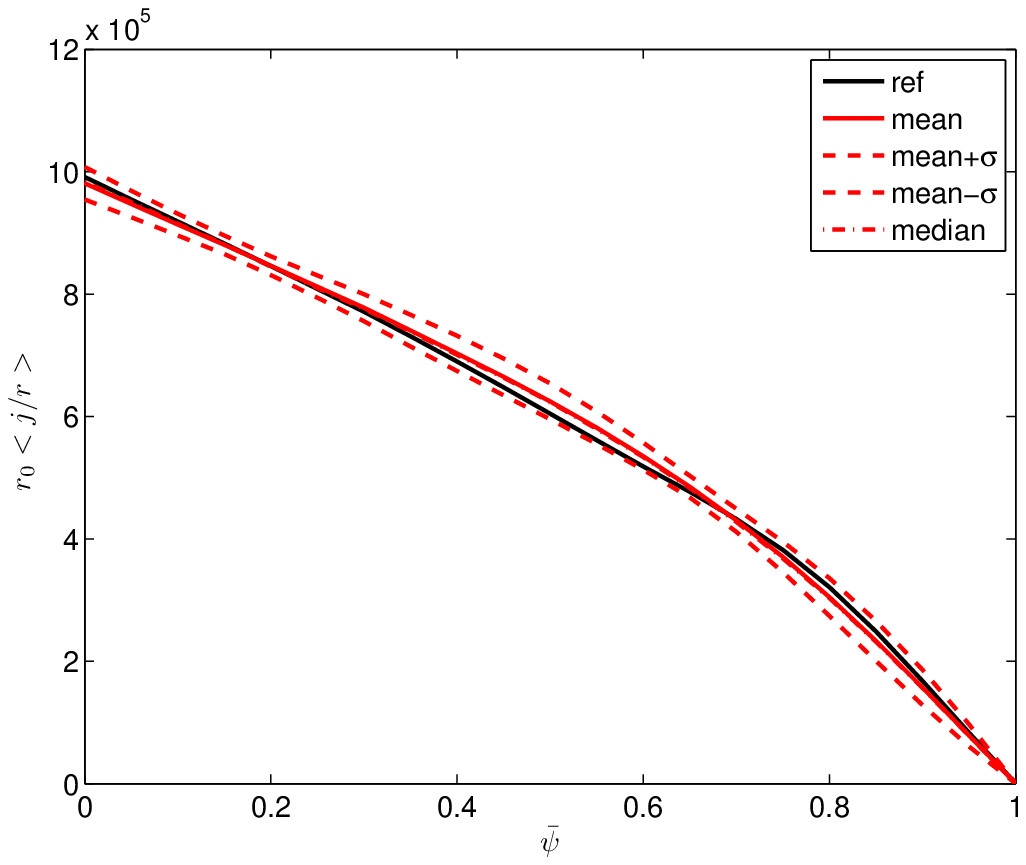}&\includegraphics[width=6cm]{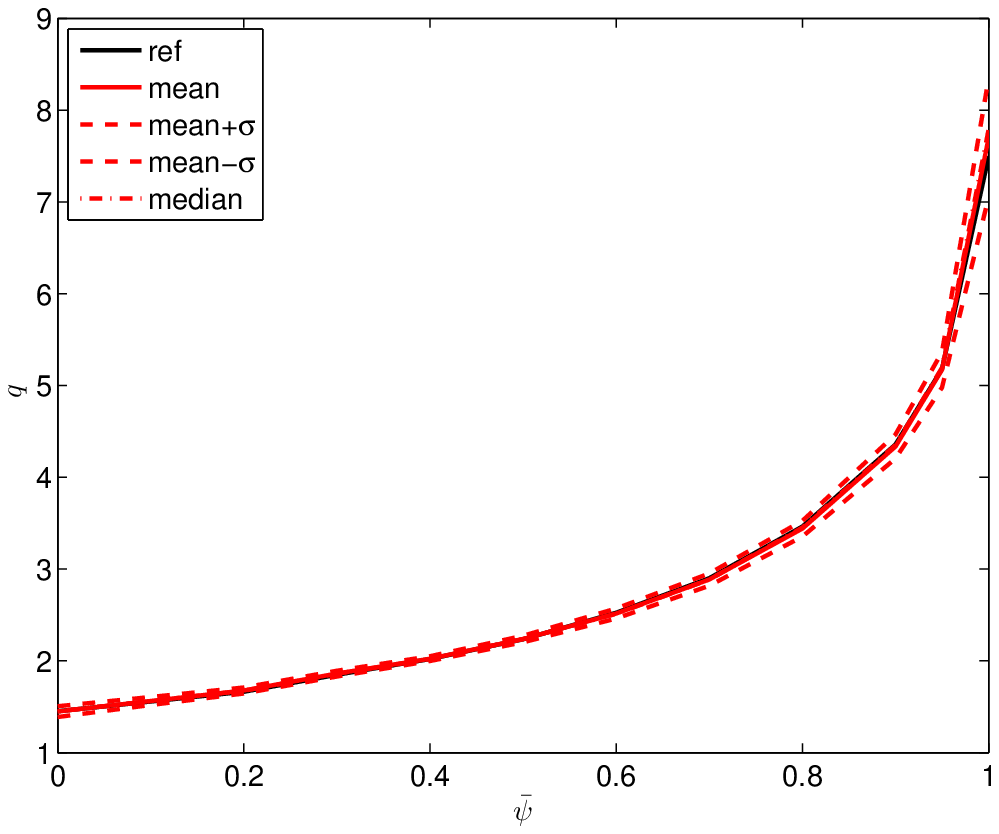} 
 \end{tabular}
 \end{center}
 \caption{
Statistical results of the identification experiments with noisy measurements 
(magnetics, interferometry and polarimetry). 
Row 1: $\eps=10^{-2}$, row 2: $\eps=10^{-1}$, row 3 $\eps=1$. Column 1: $ R_0\ds <\frac{j(r,\p)}{r}>$,
 and column 2: safety factor $q$.
For each function the reference value is given in black and the mean identified function in red. 
The mean $\pm$ standard deviation functions are shown in dashed red. \label{fig:twinmeannoisepolar2}}
 \end{figure}

The statistical results of the twin experiments are shown on 
Figs. \ref{fig:twinmeannoisepolar} and \ref{fig:twinmeannoisepolar2} for 3 different values of $\eps$. 
The use of interferometric and polarimetric measurements adds supplementary constraints on the $A$ and $B$ functions.  
The variability in the recovered functions is less important than in the case where only 
magnetics are used particularly for $\p \in [0,0.5]$. 
This is not surprising since the new measurements are internal and bring some information contained inside the 
plasma domain. Nevertheless it is not enough to perfectly reconstruct independently the $A$ and $B$ functions.
This does not prevent an excellent recovery of the averaged current density profile and of the safety factor $q$. 
This phenomenon already observed in the magnetics case is emphasized here where the variability 
of the recovered profiles has decreased.

\subsection{A real pulse}

The algorithm detailed in this paper has been implemented in a C++ software called 
Equinox developed in collaboration with the Fusion Department at
Cadarache for Tore Supra and JET. 
Equinox can be used on the one hand for precise studies in which the computing time is not a limiting factor 
and on the other hand in a real-time framework to reconstruct the successive 
plasma equilibrium configurations during a whole pulse. 
For the time being it is used on JET and ToreSupra pulses, it has also been tested on the Tokamak TCV 
and can potentially be used on any Tokamak.

During the real time analysis of a whole pulse an equilibrium is reconstructed from 
new measurements with a time step of $\Delta t=$ 100 ms. 
For each equilibrium reconstruction the number of iterations of the algorithm is set to $2$. 
This enables fast enough computations while a very good precision is achieved since the initial guess for an 
equilibrium computation at time $t$ is the equilibrium computed at time $t-\Delta t$. 
After 1 iteration a typical value for the relative residu on $\psi$ is of $10^{-2}$ and it is of $10^{-3}$ 
after 2 iterations. Table \ref{tab:realtime} gives the size of the finite elements mesh used at ToreSupra and at JET 
as well as typical computation times on a laptop computer.

\begin{table}[!h]
\begin{center}
\begin{tabular}{|l|l|l|}
\hline
  & ToreSupra & JET \\
\hline
\multicolumn{3}{|c|}{Finite element mesh}\\
\hline
Number of triangles & $1382$ & $2871$\\
\hline
Number of nodes & $722$ & $1470$ \\
\hline
\multicolumn{3}{|c|}{Computation time (1.80GHz) }\\
\hline
One equilibrium & $20$ ms & $60$ ms \\
\hline
\end{tabular}
\caption{Typical mesh size and computation time for ToreSupra and JET \label{tab:realtime} }
\end{center}
\end{table}

The choice of the regularization parameters is crucial since it determines 
the balance between the fit to the data and the regularity of the identified functions. It is also 
difficult as is shown in the preceding section. 
Ideally they should be determined for each equilibrium reconstruction. However this is not possible 
in a real-time application and the regularization parameters have to be set apriori 
to a constant value. 
From the twin experiments presented in the preceding sections it is quite clear 
that a good value for the regularization parameter $\eps$ is in the range $[10^{-2},1]$.
By trial and error on different pulses using magnetics, interferometry and polarimetry, 
it appeared that a value of $\eps=5.10^{-2}$ gave good results.

\begin{figure}[!h]
\begin{center}
  \includegraphics[width=10cm]{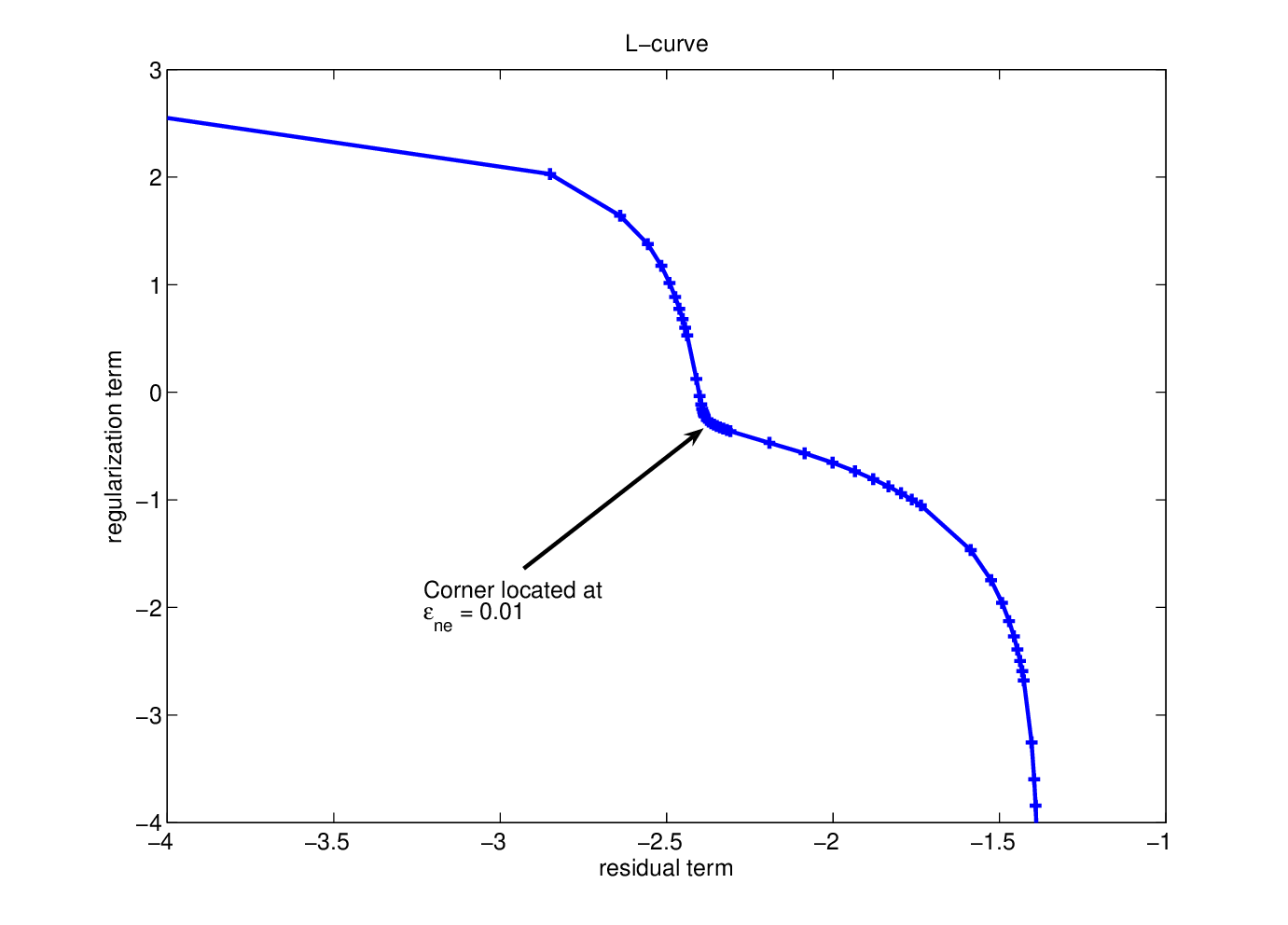} 
\end{center}
\caption{
Typical Lcurve for the determination of $\eps_{n_e}$. 
It is a plot of the parametric curve 
$x(\eps_{n_e})=log(\ds \frac{1}{2}||D^{1/2}(Bv^*(\eps_{n_e})-\gamma)||^2)$, 
$y(\eps_{n_e})=log(\frac{1}{2} (v^*(\eps_{n_e}))^T \Lambda v^*(\eps_{n_e}))$ 
where $v^*(\eps_{n_e})$ is the solution to Eq. (\ref{eqn:nenormal}). 
Hansen's algorithm \cite{Hansen:1999} locates a corner at $\eps_{n_e}=0.01$. 
\label{fig:LcurveNe}}
\end{figure}

As for the identification of functions $A$ and $B$ the choice of a good regularization 
parameter for the identification of $n_e$ is crucial. 
However in this case the L-curve method works quite well and 
it was used to determine the regularization parameters $\eps_{n_e}$ a priori 
on a number of equilibria for a few shots. The obtained values showed little 
variation and the choice of a mean value $\eps=0.01$ proved to be efficient. 
Figure \ref{fig:LcurveNe} shows an example of an L-curve computed for the identification of $n_e$.

Concerning real pulses at JET we refer to \cite{ACTI.B.Faugeras.09.1, OS.B.Faugeras.10.1} 
in which a validation of Equinox is performed using many different pulses. 
This validation includes a posteriori comparison of the position of rational q surfaces 
computed from 
Equinox and deduced from soft X-rays measurements. The validation is satisfactory and shows again that 
when solving the inverse problem the use of interferometry, polarimetry 
and even Motional Strak Effect measurements at JET improves the location of rational q surfaces.    

\begin{figure}[!h]
\begin{center}
\begin{tabular}{l}
\includegraphics[width=14cm]{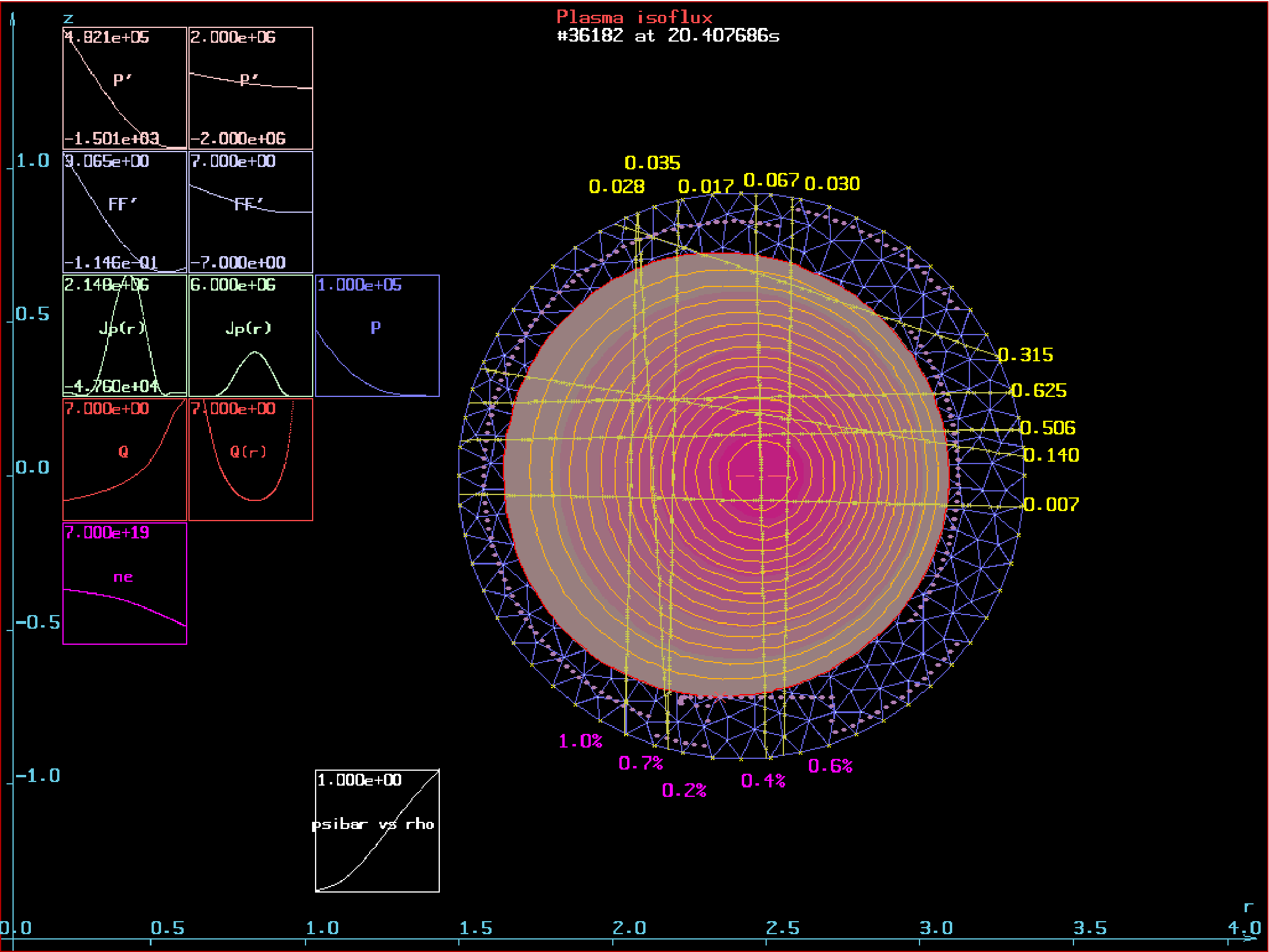} 
\end{tabular}
\end{center}
\caption{Graphical output from Equinox. Reconstructed equilibrium at time 20.408 s 
for ToreSupra pulse number 36182. Magnetic, interferometry and polarimetry measurments are used. 
See text for more details. 
\label{fig:equinox-mag-polar}}
\end{figure}

Here we only present an example of the output from Equinox on a ToreSupra pulse. 
Figure \ref{fig:equinox-mag-polar} shows the equilibrium computed at time 20.408 seconds 
for ToreSupra pulse number 36182 using magnetic measurements as well as interferometric 
and polarimetric measurements. One can observe the position of the plasma in the vacuum vessel. 
Isoflux lines are displayed from the magnetic axis to the boundary. 
The interferometry and polarimetry chords are displayed. 
For each chord the error between computed and measured interferometry is given in purple. 
These errors are about $1\%$ for the active chords. 
The polarimetry absolute errors are given in yellow. 
Different graphs are plotted on the left hand side of the display. 
On the first row the identified function $A$, and corresponding functions $p'$ and $p$. 
On the second row the identified function $B$ and corresponding function $ff'$. 
The third row gives the toroidal current density $j_{\phi}$ in the equatorial plane  
and the fourth one shows the safety factor $q$. 
Finally on the fifth row the identified $n_e$ function is plotted.

It is of importance to compute the kinetic energy poloidal $\beta_p$ parameter and the internal inductance $l_i$. 
In Equinox these equilibrium parameters are computed following the equations of Appendix \ref{appendixC}. 
For ToreSupra they are computed in the code Apolo \cite{Saint-Laurent:2001} from the Shafranov integrals and from 
the toroidal plasma flux. The agreement between the two methods is good as shown in Table \ref{tab:betali}. 
The relative errors on $\beta_p$ and $l_i$ are about $10\%$ while it is of about $1\%$ 
on the sum $\beta_p+\ds \frac{l_i}{2}$.

\begin{table}[!h]
\begin{center}
\begin{tabular}{|l|l|l|l|}
\hline
   & $\beta_p$  & $l_i$  & $\beta_p + \ds \frac{l_i}{2}$  \\
\hline
 Equinox & $0.62$ & $1.66$ & $1.45$ \\
\hline
 Apolo & $0.70$ & $1.55$ & $1.47$ \\
\hline
\end{tabular}
\caption{$\beta_p$ and $l_i$ computed by Equinox and by Apolo for ToreSupra shot 36182 at t=20.408s \label{tab:betali} }
\end{center}
\end{table}

Finally it should be noticed that at ToreSupra or JET there does not exist 
reliable enough pressure measurements to be used in an inverse equilibrium reconstuction. 
The electron pressure $p_e$ can be reasonably estimated from interferometry for the density $n_e$ 
and Thomson scattering and Electron Cyclotron Emission for the temperature $T_e$. On the contrary very large 
uncertainties on the ion quantities $n_i$ and $T_i$ make the ion pressure $p_i$ 
and thus the total pressure $p=p_e+p_i$ unusable in 
a real-time identification algorithm such as the one presented here. 
Moreover the quantity really important in order to constrain the identification of the $p'$ term 
would be the pressure gradient on which the error bars are even larger.

\section{Conclusion}
We have presented an algorithm for 
the identification of the current density profile in the Grad-Shafranov equation and 
the equilibrium reconstruction from experimental measurements in real time. 
We have shown thanks to several twin experiments that even though the unknown 
functions $A$ and $B$ (or $p'$ and $ff'$) taken separately might not be always exactly identified 
the resulting averaged current density and safety factor seem to be always well identified. 
We have also shown that the use of internal polarimetric measurements improves the quality of the 
identification but is still not enough to perfectly identify both $A$ and $B$. 
Finally we have introduced the software Equinox in which this methodology is developed. 
This work constitutes a step towards the real-time control of the safety factor 
and of the averaged current density profile in a Tokamak plasma which will be essential in nuclear fusion reactors.

\section*{Aknowledgements}
The authors are grateful to Kristoph Bosak who developed a first version of the code Equinox. 
Although it has now been thoroughly modified this version was an essential basis to start from.\\
The authors would also like to thank all colleagues from the CEA at Cadarache in France 
involved in a collaboration between the University of Nice and the CEA 
through the LRC (Laboratoire de Recherche Conventionn\'e). 
Discussions with Francois Saint-Laurent and Sylvain Bremond were particularly helpful. 
Emmanuel Joffrin initiated the real-time approach and Didier Mazon helped introducing us at JET 
where different people are also involved. In particular Luca Zabeo provided magnetic input 
data from the boundary code Xloc for Equinox and the work of Fabio Piccolo and Robert Felton 
is essential to implement Equinox on JET real-time system.

\clearpage
\newpage


\appendix

\section{Average over magnetic surfaces}
\label{appendixA}
The method of averaging over the magnetic surfaces is detailed in \cite{Blum:1989} (p 242). 
The average $<A>$ of an arbitrary quantity $A$ on a magnetic surface $S$ is defined as 
$$
<A>=\ds \frac{\partial}{\partial V} \int_V A dV
$$  
where $V$ is the volume inside $S$. This notion of average has the following property:
$$
<A>=\ds \frac{\ds \int_{C_\p} \ds \frac{Adl}{B_p}}{\ds \int_{C_\p} \ds \frac{dl}{B_p}}
$$
where $C_\p$ is a closed contour ${\p=cte} \in (0,1)$ and $B_p=\ds \frac{1}{r}||\nabla \psi||$.

\section{Safety factor $q$}
\label{appendixB}
The safety factor is so called because of the role it plays 
in determining stability (\cite{Wesson:2004} p 111). 
It can be seen as the ratio of the variation of the toroidal 
angle needed for one magnetic field line to perform one poloidal turn.
$$
q=\frac{\Delta \phi}{2\pi}
$$
Since $q$ is the same for all magnetic field lines on a magnetic surface it is a function of $\psi$ 
(or $\p$). The expression of $q$ used for computations is the following
\begin{equation*}
q(\p)=\frac{1}{2\pi} \int_{C_\p} \frac{B_{\phi}}{rB_p}dl
\end{equation*}
where $C_\p$ is a closed contour ${\p=cte} \in (0,1)$, $B_\phi=\ds \frac{f}{r}$
and
$$
f(\psi)= \sqrt{ (B_0 R_0)^2 + \ds \int_{\psi_b}^{\psi} (f^2)'(y)dy }
$$


\section{Poloidal $\beta_p$ and Internal inductance $l_i$}
\label{appendixC}
The full 3D plasma domain is denoted by $D$. 
The plasma domain in the poloidal section by $\Omega_p$ and its boundary $\partial \Omega_p = \Gamma_p$. Let us define $R_{g}=\ds \frac{1}{2}(R_{left}+R_{right})$

\paragraph{Surface and perimeter of a poloidal section}
Let us define $S_p=\int_{\Omega_p} ds $ and $ L_p=\int_{\Gamma_p} dl $. 
For a circular plasma of radius $a$: $L_p=2\pi a$, $S_p=\pi a^2 $ and $S_p=\ds \frac{L_p^2}{4\pi}$.
Even for non-circular plasma the following quantity is used:
\begin{equation}\label{eqn:surfaceapprox}
\hat{S}_p=\frac{L_p^2}{4\pi} 
\end{equation}

\paragraph{Plasma volume}
\begin{equation}
V_p=\int_D dv = \int_0^{2\pi} \int_{\Omega_p} r d\phi ds = 2\pi \int_{\Omega_p}rds 
\end{equation}
The following approximation can be used:
\begin{equation}\label{eqn:volumeapprox}
\hat{V}_p=2 \pi R_{g} \hat{S}_p 
\end{equation}

\paragraph{Poloidal $\beta_p$}
The ratio $\beta=\ds \frac{p}{B^2/2\mu_0}$ represents the efficiency of the confinement of the plasma pressure
by the magnetic field. The poloidal beta is defined as the ratio of the mean kinetic pressure of the plasma to 
its magnetic pressure (\cite{Wesson:2004} p 116): 
\begin{equation}
\beta_p= \frac{\bar{p}}{B_{pa}^2/2\mu_0}
\end{equation} 
where
\begin{equation}
\bar{p}=\frac{\int_{D} p dv}{\int_{D}dv}=\frac{\int_{\Omega_p} p rds}{\int_{\Omega_p} rds}
\end{equation}
and
\begin{equation}\label{eqn:Bpa}
B_{pa}=\frac{\int_{\Gamma_p} B_p dl}{\int_{\Gamma_p} dl} = \frac{\mu_0 I_p}{L_p}
\end{equation}
Let us define the internal kinetic energy
$$
W=\frac{3}{2}\int_Dpdv
$$ 
We have 
$$
W=\frac{3}{2}\bar{p}V_p=\frac{3}{2}\frac{B_{pa}^2}{2\mu_0} V_p \beta_p
$$
and from Eq. (\ref{eqn:Bpa}), (\ref{eqn:volumeapprox}) and (\ref{eqn:surfaceapprox}) 
follows that (\cite{Wesson:2004} p 504)
$$
W = \ds \frac{3}{8}\mu_0 R_g I_p^2 \beta_b
$$
Then $\beta_p$ can be approximated by 
\begin{equation}
\label{betaapprox}
\beta_p=\ds \frac{\ds \frac{3}{2}\bar{p}V_p}{\ds \frac{3}{8}\mu_0 R_g I_p^2}
\end{equation}
which the default $\beta_p$ computed by Equinox.

\subsection{Internal inductance $l_i$}
The internal inductance $l_i$ of the plasma characterizes the current density profile 
(\cite{Wesson:2004} p 120, \cite{Blum:1989} p 44):
\begin{equation}
l_i=\frac{\bar{B}_p^2}{B_{pa}^2}
\end{equation}
where
\begin{equation*}
\bar{B}_p^2=\frac{\int_D B_p^2 dv}{\int_D dv}
\end{equation*}
In Equinox the computation of $l_i$ is done as follows:
\begin{equation*}
l_i=\frac{\bar{B}_p^2 V_p}{B_{pa}^2 V_p}
\end{equation*}
Using Eq. (\ref{eqn:Bpa}), Eq. (\ref{eqn:volumeapprox}) and Eq. (\ref{eqn:surfaceapprox}) 
leads to
\begin{equation}
l_i = \ds \frac{\bar{B}_p^2 V_p}{\ds \frac{\mu_0^2}{2}R_g I_p^2}
\end{equation} 
which is the default computation of $l_i$ in Equinox.

\clearpage
\newpage
\bibliographystyle{elsarticle-num}
\bibliography{biblio-fusion}

\end{document}